\documentclass[twoside, 12pt]{article}
\usepackage{amssymb,mathrsfs,psfrag,graphicx}
\usepackage{hyperref}
\def\v{\varepsilon}

\def\x{\xi}
\def\t{\theta}
\def\T{\Theta}

\def\mb{\mathbf}
\def\a{\alpha}

\def\g{\gamma}
\def\d{\delta}

\def\f{\frac}
\def\p{\phi}

\def\z{\zeta}

\def\di{\displaystyle}
\def\i{\infty}

\begin{document}

\title{\bf  Hydrodynamic Limit of the Boltzmann Equation with Contact Discontinuities}
\vskip 0.5cm
\author{{\bf Feimin Huang$^*$, Yi Wang$^*$, Tong Yang$^{\dag}$}\\
$^*$Institute of Applied Mathematics,\\ AMSS, Academia Sinica, China\\
$^{\dag}$ Department of Mathematics, \\City University of HongKong,
HongKong
 }

\date{ }
\maketitle

\begin{abstract}
The hydrodynamic limit  for the  Boltzmann equation
is studied in the case when the limit system, that is, the
 system of Euler equations
contains contact discontinuities. When suitable
initial data is chosen
to avoid the initial layer, we prove that there
exists a unique solution to the Boltzmann equation globally in time for any
given Knudsen number. And this family of solutions
converge to the local Maxwellian defined by the
contact discontinuity  of the
Euler equations uniformly away from the discontinuity as the Knudsen
number $\v$ tends to zero.
The proof is obtained by an appropriately chosen scaling  and the
energy method through  the micro-macro decomposition.
\end{abstract}

\section{Introduction}
\renewcommand{\theequation}{\arabic{section}.\arabic{equation}}
\setcounter{equation}{0}
Consider the Botlzmann equation with
 slab symmetry
\begin{equation}
f_t + \xi_1f_x= \f{1}{\v}Q(f,f),~(f,x,t,\xi)\in {\mb R}\times {\mb
R}\times{\mb R}^+
 \times {\mb R}^3, \label{(1.1)}
\end{equation}
where $\xi=(\xi_1,\xi_2,\xi_3)\in {\mb R}^3$, $f(x,t,\xi)$ is the
 density distribution function of the particles at time $t$ and
space $x$ with velocity $\xi$,  and
$\v>0$ is the Knudsen number which is proportional to the mean free
path.

The equation (\ref{(1.1)}) was established by Boltzmann \cite{Boltzmann}
in 1872 to describe the motion of  rarefied gases and it is
a fundamental equation in statistics
physics. For monatomic gas, the rotational invariance of the
particles leads to the following  bilinear
form for the collision operator
$$
 Q(f,g)(\xi) \equiv \f{1}{2}
\int_{{\mb R}^3}\!\!\int_{{\mb S}_+^2} \Big(f(\xi')g(\xi_*')+
f(\xi_*')g(\xi')-f(\xi) g(\xi_*) - f(\xi_*)g(\xi) \Big)
B(|\xi-\xi_*|, \hat\t)
 \; d \xi_* d \Omega,
$$
where $\xi',\xi_*'$ are the velocities after an elastic collision of
two
particles with velocities  $\xi,\xi_*$ before the collision. Here,
$\hat\t$ is the angle between the relative velocity $\xi-\xi_*$
and the unit vector $\Omega$ in ${\mb S}^2_+=\{\Omega\in {\mb
S}^2:\ (\xi-\xi_*)\cdot \Omega\geq 0\}$. The conservation of
momentum and energy gives the following relation
between the velocities before and after collision:
$$
\left\{
\begin{array}{l}
 \xi'= \xi -[(\xi-\xi_*)\cdot \Omega] \; \Omega, \\[3mm]
 \xi_*'= \xi_* + [(\xi-\xi_*)\cdot \Omega] \; \Omega.
\end{array}
\right.
$$

In this paper, we consider the Boltzmann equation for the two basic models,
that is, the hard sphere model and the hard potential including Maxwellian
molecules under the assumption of angular
cut-off. That is, we assume that the collision kernel
$B(|\xi-\xi_*|,\hat\t)$ takes one of the following two forms,
$$
B(|\xi-\xi_*|,\hat\t)=|(\xi-\xi_*, \Omega)|,
$$
and
$$
B(|\xi-\xi_*|,\hat\t)=|\xi-\xi_*|^{\f{n-5}{n-1}}b(\hat\t),\quad
b(\hat\t)\in L^1([0, \pi]),~n\ge 5.
$$
 Here, $n$ is the index in the inverse power potential which is
proportional to $r^{1-n}$ with $r$ being the distance between two
particles.

Formally,  when the Knudsen number $\v$ tends to zero,
the limit  of the Boltzmann equation (\ref{(1.1)}) is
 the classical system of Euler equations
\begin{equation}
\left\{
\begin{array}{l}
\di\rho_t+(\rho u_1)_x=0,\\
\di(\rho u_1)_t+(\rho u_1^2+p)_x=0,\\
\di(\rho u_i)_t+(\rho u_1u_i)_x=0,~i=2,3,\\
\di[\rho(E+\f{|u|^2}{2})]_t+[\rho u_1(E+\f{|u|^2}{2})+pu_1]_x=0,
\end{array}
\right. \label{(1.2)}
\end{equation}
where
\begin{equation}
\left\{
\begin{array}{l}
\di\rho(x,t)=\int_{\mb{R}^3}\varphi_0(\xi)f(x,t,\xi)d\xi,\\
\di\rho u_i(x,t)=\int_{\mb{R}^3}\varphi_i(\xi)f(x,t,\xi)d\xi,~i=1,2,3,\\
\di\rho(E+\f{|u|^2}{2})(x,t)=\int_{\mb{R}^3}\varphi_4(\xi)f(x,t,\xi)d\xi.
\end{array}
\right. \label{(1.3)}
\end{equation}
Here, $\rho$ is the density, $u=(u_1,u_2,u_3)$ is the macroscopic
velocity, $E$ is the internal energy and the pressure $p=R\rho\t$
with $R$ being the gas constant. The temperature $\t$ is related to the
internal energy by $E=\f{3}{2}R\t$, and
$\varphi_i(\xi)(i=0,1,2,3,4)$ are
the collision invariants given by
$$
\left\{
\begin{array}{l}
 \varphi_0(\xi) = 1, \\
 \varphi_i(\xi) = \xi_i \ \ {\textrm {for} }\ \  i=1,2,3, \\
 \varphi_4(\xi) = \f{1}{2} |\xi|^2,
\end{array}
\right.
$$
that satisfy
$$
\int_{{\mb R}^3} \varphi_i(\xi) Q(h,g) d \xi =0,\quad {\textrm {for}
} \ \  i=0,1,2,3,4.
$$

How to justify the above limit, that is,  the Euler equation
(\ref{(1.2)}) from Boltzmann equation (\ref{(1.1)}) when Knudsen
number tends to zero is an open problem going way back to
the time of Maxwell. For this,  Hilbert  introduced
the famous Hilbert expansion to show  formally that  the first order
approximation of the Boltzmann equation gives  the Euler equations.
On the other hand, it is important to verify this limit
process rigorously in mathematics. For the case when the Euler equation has
smooth solutions,
the zero Knudsen number limit of the Boltzmann equation has been studied
even in the
case with an initial layer, cf.  Asona-Ukai\cite{Asona-Ukai},
Caflish\cite{Caflish}, Lachowicz\cite{Lachowicz} and
Nishida\cite{Nishida} etc. However, as is well-known,
solutions of the Euler equation (\ref{(1.2)}) in general develop  singularities, such as shock waves and contact
discontinuities. Therefore,  how to
verify the hydrodynamic limit from Boltzmann equation to the Euler
equations with basic wave patterns is an natural problem.
In this direction,  Yu \cite{Yu} showed that
when the solution of the Euler equation (\ref{(1.2)}) contains non-interacting
shocks,
there exists a sequence of solutions to
 the Boltzmann equation that converge to the local Maxwellian defined by the
 solution of the Euler equation (\ref{(1.2)}) uniformly
away from the shock. In this work, the inner and outer expansions
developed by Goodman-Xin \cite{Goodman-Xin} for conservation laws and
the Hilbert expansion were crucially used.

The main purpose of this paper is to study the hydrodynamic
limit of the Boltzmann equation when the
corresponding Euler equation contains contact discontinuities. More
precisely, given a  solution of the Euler
equation (\ref{(1.2)}) with contact discontinuities, we will show that there exists a  family of solutions to
the Boltzmann equation that  converge to a local Maxwellian defined by
the Euler solution
 uniformly away from the contact
discontinuity as $\v\rightarrow 0$. Moreover, a uniform convergence
rate in $\v$ is also given. The proof is obtained by a scaling
transformation of the independent variables and the perturbation
together with  the energy method introduced
 by Liu-Yang-Yu\cite{Liu-Yang-Yu}.

For later use, we now  briefly introduce the micro-macro decomposition around the
local Maxwellian defined by the solution to the Boltzmann equation,
cf.
\cite{Liu-Yang-Yu}. For a solution $f(x,t,\xi)$ of the Boltzmann
equation (\ref{(1.1)}), we decompose it into
$$
f(x,t,\xi)=\mb{M}(x,t,\xi)+\mb{G}(x,t,\xi),
$$
where the local Maxwellian
$\mb{M}(x,t,\xi)=\mb{M}_{[\rho,u,\t]}(\xi)$ represents the
macroscopic (fluid) component of the solution, which is naturally
defined by the five conserved quantities, i.e., the mass
density $\rho(x,t)$, the momentum $\rho u(x,t)$, and the total
energy $\rho(E+\f{1}{2}|u|^2)(x,t)$ in (\ref{(1.3)}), through
\begin{equation}
\mb{M}=\mb{M}_{[\rho,u,\t]} (x,t,\xi) = \f{\rho(x,t)}{\sqrt{ (2 \pi
R \t(x,t))^3}} e^{-\f{|\xi-u(x,t)|^2}{2R\t(x,t)}}. \label{(1.4)}
\end{equation}
And $\mb{G}(x,t,\xi)$ being the difference between the solution and
the above local Maxwellian represents the microscopic (non-fluid) component.

For convenience, we denote the inner product of $h$ and $g$ in
$L^2_{\xi}({\mb R}^3)$ with respect to a given Maxwellian
$\tilde{\mb{M}}$ by:
$$
 \langle h,g\rangle_{\tilde{\mb{M}}}\equiv \int_{{\mb R}^3}
 \f{1}{\tilde{\mb{M}}}h(\xi)g(\xi)d \xi.
$$
 If $\tilde{\mb{M}}$ is the local
Maxwellian $\mb{M}$ defined in (\ref{(1.4)}), with respect to the
corresponding inner product, the macroscopic space is spanned by the
following five pairwise orthogonal base
$$
\left\{
\begin{array}{l}
 \chi_0(\xi) \equiv {\di\f1{\sqrt{\rho}}\mb{M}}, \\[2mm]
 \chi_i(\xi) \equiv {\di\f{\xi_i-u_i}{\sqrt{R\t\rho}}\mb{M}} \ \ {\textrm {for} }\ \  i=1,2,3, \\[2mm]
 \chi_4(\xi) \equiv
 {\di\f{1}{\sqrt{6\rho}}(\f{|\xi-u|^2}{R\t}-3)\mb{M}},\\
 \langle\chi_i,\chi_j\rangle=\delta_{ij}, ~i,j=0,1,2,3,4.
 \end{array}
\right.
$$
In the following, if $\tilde{\mb{M}}$ is the local Maxwellian $\mb{M}$, we
just use the simplified notation $\langle\cdot,\cdot\rangle$ to denote
 the
inner product $\langle\cdot,\cdot\rangle_{\mb{M}}$.
We can now define the macroscopic projection $\mb{P}_0$ and
microscopic projection $\mb{P}_1$ as follows
\begin{equation}
\left\{
\begin{array}{l}
 \mb{P}_0h = {\di\sum_{j=0}^4\langle h,\chi_j\rangle\chi_j,} \\
 \mb{P}_1h= h-\mb{P}_0h.
 \end{array}
\right.\label{(1.5)}
\end{equation}
The projections $\mb{P}_0$ and $\mb{P}_1$ are orthogonal and satisfy
$$
\mb{P}_0\mb{P}_0=\mb{P}_0, \mb{P}_1\mb{P}_1=\mb{P}_1,
\mb{P}_0\mb{P}_1=\mb{P}_1\mb{P}_0=0.
$$
We remark that a function $h(\xi)$ is called microscopic or
non-fluid if
$$
\int h(\xi)\varphi_i(\xi)d\xi=0,~i=0,1,2,3,4,
$$
where $\varphi_i(\xi)$ is the collision invariants.

 Under the above micro-macro decomposition, the solution $f(x,t,\xi)$ of the Boltzmann
equation (\ref{(1.1)}) satisfies
$$
\mb{P}_0f=\mb{M},~\mb{P}_1f=\mb{G},
$$
and the Boltzmann equation (\ref{(1.1)}) becomes
\begin{equation}
(\mb{M}+\mb{G})_t+\xi_1(\mb{M}+\mb{G})_x
=\f{1}{\v}[2Q(\mb{M},\mb{G})+Q(\mb{G},\mb{G})]. \label{(1.6)}
\end{equation}
If we multiply the equation (\ref{(1.6)}) by the collision invariants
$\varphi_i(\xi)(i=0,1,2,3,4)$ and integrate the resulting equations
with respect to $\xi$ over ${\mb R}^3$, then we can get
the following fluid-type system for the fluid components:
\begin{equation}
\left\{
\begin{array}{lll}
\di \rho_{t}+(\rho u_1)_x=0, \\
\di (\rho u_1)_t+(\rho u_1^2
+p)_x=-\int\xi_1^2\mb{G}_xd\xi,  \\
\di (\rho u_i)_t+(\rho u_1u_i)_x=-\int\xi_1\xi_i\mb{G}_xd\xi,~ i=2,3,\\
\di [\rho(E+\f{|u|^2}{2})]_t+[\rho
u_1(E+\f{|u|^2}{2})+pu_1]_x=-\int\f12\xi_1|\xi|^2\mb{G}_xd\xi.
\end{array}
\right. \label{(1.7)}
\end{equation}

 Note that the above fluid-type system is not
closed and we need  one more  equation for the non-fluid
component ${\mb{G}}$ which can be obtained
by applying the projection operator $\mb{P}_1$
to the equation (\ref{(1.6)}):
\begin{equation}
\mb{G}_t+\mb{P}_1(\xi_1\mb{M}_x)+\mb{P}_1(\xi_1\mb{G}_x)
=\f{1}{\v}\left[\mb{L}_\mb{M}\mb{G}+Q(\mb{G}, \mb{G})\right].
\label{(1.8)}
\end{equation}
Here $\mb{L}_\mb{M}$ is the linearized collision operator of
$Q(f,f)$
 with respect to the local Maxwellian $\mb{M}$:
$$
\mb{L}_\mb{M} h=2Q(\mb{M}, h)=Q(\mb{M}, h)+ Q(h,\mb{M}).
$$
And the null space $\mathcal{N}$ of $\mb{L}_\mb{M}$ is spanned by
the macroscopic variables:
$$
\chi_j(\xi), ~j=0,1,2,3,4.
$$
Furthermore, there exists a positive constant
$\sigma_0(\rho,u,\t)>0$ such that for any function $h(\xi)\in
\mathcal{N}^\bot$, cf. \cite{Grad},
$$
<h,\mb{L}_\mb{M}h>\le -\sigma_0<\nu(|\xi|)h,h>,
$$
where $\nu(|\xi|)$ is the collision frequency. For the hard sphere
and the hard potential with angular cut-off, the collision
frequency $\nu(|\xi|)$ has the following property
$$
0<\nu_0<\nu(|\xi|)\le c(1+|\xi|)^{\kappa},
$$
for some positive constants $\nu_0, c$ and $0\le\kappa\le 1$.

Consequently, the linearized collision operator $\mb{L}_\mb{M}$ is a
dissipative operator on $L^2({\mb R}^3)$, and its inverse
$\mb{L}_\mb{M}^{-1}$  exists and is  a bounded
operator in $L^2({\mb R}^3)$.

It follows from (\ref{(1.8)}) that
\begin{equation}
\mb{G}=\v \mb{L}_\mb{M}^{-1}[\mb{P}_1(\xi_1\mb{M}_x)] +\T,
\label{(1.9)}
\end{equation}
with
\begin{equation}
\T=\mb{L}_\mb{M}^{-1}[\v(\mb{G}_t+\mb
{P}_1(\xi_1\mb{G}_x))-Q(\mb{G}, \mb{G})].\label{(1.10)}
\end{equation}

Plugging the equation (\ref{(1.9)}) into (\ref{(1.7)}) gives
\begin{equation}
\left\{
\begin{array}{l}
\di \rho_{t}+(\rho u_1)_x=0,\\
\di (\rho u_1)_t+(\rho u_1^2 +p)_x=\f{4\v}{3}(\mu(\t)
u_{1x})_x-\int\xi_1^2\T_xd\xi,  \\
\di (\rho u_i)_t+(\rho u_1u_i)_x=\v(\mu(\t)
u_{ix})_x-\int\xi_1\xi_i\T_xd\xi,~ i=2,3,\\
\di [\rho(E+\f{|u|^2}{2})]_t+[\rho
u_1(E+\f{|u|^2}{2})+pu_1]_x=\v(\lambda(\t)\t_x)_x+\f{4\v}{3}(\mu(\t)u_1u_{1x})_x\\
\di +\v\sum_{i=2}^3(\mu(\t)u_iu_{ix})_x
-\int\f12\xi_1|\xi|^2\T_xd\xi,
\end{array}
\right. \label{(1.11)}
\end{equation}
where the viscosity coefficient $\mu(\t)>0$ and the heat
conductivity coefficient $\lambda(\t)>0$ are smooth functions of the
temperature $\t$, and we normalize the gas constant $R$ to be
$\f{2}{3}$ so that $E=\t$ and $p=\f23\rho\t$. The explicit formula
of $\mu(\t)$ and $\lambda(\t)$ can be found for example in \cite{Yang-Zhao}, we
omit it here for brevity.

Since the problem considered in
this paper is  one dimensional in the space variable $x\in {\bf R}$, in the
macroscopic level, it is more convenient to rewrite the equation
(\ref{(1.1)}) and the system (\ref{(1.2)}) in the {\it Lagrangian}
coordinates as in the study of conservation laws. That is, set the
coordinate transformation:
$$
 x\Rightarrow \int_0^x \rho(y,t)dy, \qquad t\Rightarrow t.
$$
We will still denote the {\it Lagrangian} coordinates by $(x, t)$
for simplicity of notation. Then (\ref{(1.1)}) and (\ref{(1.2)}) in
the Lagrangian coordinates become, respectively,
\begin{equation}
f_t-\f{u_1}{v}f_x+\f{\xi_1}{v}f_x=\f{1}{\v}Q(f,f),\label{(1.12)}
\end{equation}
and
\begin{equation}
\left\{
\begin{array}{llll}
\di v_{t}-u_{1x}=0,\\
\di u_{1t}+p_x=0,\\
\di u_{it}=0, ~i=2,3,\\
\di (\t+\f{|u|^{2}}{2}\bigr)_t+ (pu_1)_x=0.\\
\end{array}
\right.\label{(1.13)}
\end{equation}
Also, (\ref{(1.7)})-(\ref{(1.11)}) take the form
\begin{equation}
\left\{
\begin{array}{llll}
\di v_t-u_{1x}=0,\\
\di u_{1t}+p_x=-\int\xi_1^2\mb{G}_xd\xi,\\
\di u_{it}=-\int\xi_1\xi_i\mb{G}_xd\xi,
~i=2,3,\\
\di\bigl(\t+\f{|u|^{2}}{2}\bigr)_{t}+ (pu_1)_x=-\int\f12\xi_1|\xi|^2\mb{G}_xd\xi,\\
\end{array}
\right.\label{(1.14)}
\end{equation}
\begin{equation}
\mb{G}_t-\f{u_1}{v}\mb{G}_x+\f{1}{v}\mb{P}_1(\xi_1\mb{M}_x)+\f{1}{v}\mb{P}_1(\xi_1\mb{G}_x)=\f{1}{\v}(\mb{L}_\mb{M}\mb{G}+Q(\mb{G},\mb{G})),\label{(1.15)}
\end{equation}
with
\begin{equation}
\mb{G}=\v \mb{L}^{-1}_\mb{M}(\f{1}{v} \mb{P}_1(\xi_1
\mb{M}_x))+\T_1,\label{(1.16)}
\end{equation}
\begin{equation}
\T_1=\mb{L}_\mb{M}^{-1}[\v(\mb{G}_t-\f{u_1}v\mb{G}_x+\f{1}{v}\mb{P}_1(\xi_1\mb{G}_x))-Q(\mb{G},\mb{G})].\label{(1.17)}
\end{equation}
and
\begin{equation}
\left\{
\begin{array}{llll}
\di v_t-u_{1x}=0,\\
\di u_{1t}+p_x=\f{4\v}{3}(\f {\mu(\t)}vu_{1x})_{x}-\int\xi_1^2\T_{1x}d\xi,\\
\di u_{it}=\v(\f{\mu(\t)}{v}u_{ix})_x-\int\xi_1\xi_i\T_{1x}d\xi,
~i=2,3,\\
\di\bigl(\t+\f{|u|^{2}}{2}\bigr)_{t}+
(pu_1)_{x}=\v(\f{\lambda(\t)}{v}\t_x)_x+\f{4\v}{3}(\f{\mu(\t)}{v}u_1u_{1x})_x\\
\di\qquad+\v\sum_{i=2}^3(\f{\mu(\t)}{v} u_iu_{ix})_x
-\int\f12\xi_1|\xi|^2\T_{1x}d\xi.
\end{array}
\right. \label{(1.18)}
\end{equation}

In the following sections, we will apply some scaling and energy method
to these equations.

\section{Main result}
\setcounter{equation}{0}

We will state the main result in this section. For this,
 we firstly recall the construction of the
contact wave $(\bar v,\bar u,\bar\t)(x,t)$ for the Boltzmann
equation in \cite{Huang-Xin-Yang}. Consider the Euler system
(\ref{(1.13)}) with a Riemann initial data
\begin{equation}
(v,u,\t)(x,0)=\left\{
\begin{array}{l}
(v_-,0,\t_-),~~~x<0,\\
(v_+,0,\t_+),~~~x>0,
\end{array}
\right. \label{(2.1)}
\end{equation}
where $v_\pm,\t_\pm$ are positive constant. It is
well-known (cf. \cite{Smoller}) that the Riemann problem
(\ref{(1.13)}), (\ref{(2.1)}) admits a contact discontinuity
solution
\begin{equation}
(\bar V,\bar U,\bar\T)(x,t)=\left\{
\begin{array}{l}
(v_-,0,\t_-),~~~x<0,\\
(v_+,0,\t_+),~~~x>0,
\end{array}
\right. \label{(2.2)}
\end{equation}
provided that
\begin{equation}
p_-:=\f{R\t_-}{v_-}=p_+:=\f{R\t_+}{v_+}. \label{(2.3)}
\end{equation}
Motivated by (\ref{(2.2)}) and (\ref{(2.3)}), we expect that for the
contact wave $(\bar v,\bar u,\bar\t)(x,t)$,
$$
\bar p=\f{R\bar\t}{\bar v}\approx p_+,~~~|\bar u|^2\ll1.
$$
Then the leading order of the energy equation $(\ref{(1.18)})_4$ is
\begin{equation}
\t_t+p_+u_{1x}=\v(\f{\lambda(\t)\t_x}{v})_x. \label{(2.4)}
\end{equation}
By using the mass equation $(\ref{(1.18)})_1$ and
$v\approx\f{R\t}{p_+}$, we obtain the following nonlinear diffusion
equation
\begin{equation}
\t_t=\v(a(\t)\t_x)_x,~~~a(\t)=\f{9 p_+\lambda(\t)}{10\t}.
\label{(2.5)}
\end{equation}
From \cite{Atkinson-Peletier}, \cite{Duyn-Peletier}, we know that
the nonlinear
diffusion equation (\ref{(2.5)}) admits a unique self-similar
solution $\hat{\T}(\eta),~\eta=\f{x}{\sqrt{\v(1+t)}}$ with the
following boundary conditions
$$
\hat{\T}(-\i,t)=\t_-,~~\hat{\T}(+\i,t)=\t_+.
$$
Let $\delta=|\t_+-\t_-|$. $\hat{\T}(x,t)$ has the
property
\begin{equation}
\hat{\T}_x(x,t)=\f{O(1)\delta}{\sqrt{\v(1+t)}}e^{-\f{cx^2}{\v(1+t)}},~~~~~~{\rm
as}~~~x\rightarrow\pm\i, \label{(2.6)}
\end{equation}
with some positive constant $c$  depending only on  $\t_{\pm}$.

Now  the contact wave $(\bar v,\bar u,\bar\t)(x,t)$ can be defined by
\begin{equation}
\bar v=\f{2}{3p_+}\hat{\T},~~~\bar u_1=\f{2\v
a(\hat{\T})}{3p_+}\hat{\T}_x,~~~\bar
u_i=0,(i=2,3),~~~\bar\t=\hat{\T}-\f{|\bar u|^2}{2}. \label{(2.7)}
\end{equation}
Note that $(\bar v,\bar u,\bar \t)(x,t)$ satisfies the
following system
\begin{equation}
\left\{\begin{array}{llll}
\di \bar v_t-\bar u_{1x}=0,\\
\di \bar u_{1t}+\bar p_x=\f{4\v}{3}(\f{\mu(\bar{\t})}{\bar v}\bar u_{1x})_x+R_{1x},\\
\di \bar u_{it}=\v(\f{\mu(\bar{\t)}}{\bar{v}}\bar{u}_{ix})_x,
i=2,3,\\
\di\bigl(\bar{\t}+\f{|\bar{u}|^{2}}{2}\bigr)_{t}+
(\bar{p}\bar{u}_1)_{x}=\v(\f{\lambda(\bar{\t})}{\bar{v}}\bar{\t}_x)_x+\f{4\v}{3}(\f{\mu(\bar{\t})}{\bar{v}}
\bar{u}_1\bar{u}_{1x})_x\\
\quad\di+\v(\sum_{i=2}^3\f{\mu(\bar{\t})}{\bar{v}}\bar{u}_i\bar{u}_{ix})_x+R_{2x},
\end{array}
\right. \label{(2.8)}
\end{equation}
where
\begin{equation}
R_1=\f{2\v}{3p_+}a(\hat\T)\hat\T_t+
(\bar{p}-p_+)-\f{4\v\mu(\bar{\t})}{3\bar{v}}\bar{u}_{1x}=\di
O(1)\delta\v(1+t)^{-1}e^{-\f{cx^2}{\v(1+t)}},\label{(2.9)}
\end{equation}
\begin{equation}
\begin{array}{ll}
 R_2&\di=
\f{\v}{\bar{v}}(\lambda(\hat\T)\hat\T_x-\lambda(\bar{\t})\bar{\t}_x)
+(\bar{p}-p_+)\bar{u}_1-\f{4\v\mu(\bar{\t})}{3\bar{v}}
\bar{u}_1\bar{u}_{1x}\\
&\di=O(1)\delta\v^{3/2}(1+t)^{-3/2}e^{-\f{cx^2}{\v(1+t)}},
\end{array}
\label{(2.10)}
\end{equation}
with some positive constant $c>0$  depending only on $\t_\pm$.

From (\ref{(2.6)}), we have
\begin{equation}
\left\{
\begin{array}{l}
|\hat\T-\t_-|= O(1)\delta e^{-\f{cx^2}{2\v(1+t)}},~~~~~{\rm if}~x<0,\\
|\hat\T-\t_+|= O(1)\delta e^{-\f{cx^2}{2\v(1+t)}},~~~~~{\rm if}~x>0.
\end{array}
\right. \label{(2.11)}
\end{equation}
Therefore,
\begin{equation}
\left\{
\begin{array}{l}
|(\bar v,\bar u,\bar\t)(x,t)-(v_-,0,\t_-)|= O(1)\delta e^{-\f{cx^2}{2\v(1+t)}},~~~~~{\rm if}~x<0,\\
|(\bar v,\bar u,\bar\t)(x,t)-(v_+,0,\t_+)|= O(1)\delta
e^{-\f{cx^2}{2\v(1+t)}},~~~~~{\rm if}~x>0.
\end{array}
\right. \label{(2.12)}
\end{equation}

We are now ready to  state the main result as follows.

\

{\it \noindent{\bf Theorem 2.1} Given a contact discontinuity
solution $(\bar V,\bar U,\bar\T)(x,t)$ of the Euler system
(\ref{(1.13)}), there exists  small positive constants $\delta_0$,
$\v_0$ and a global Maxwellian $\mb{M}_*=\mb{M}_{[\rho_*,u_*,\t_*]}$,
such that if $\delta\leq \delta_0$, $\v\leq \v_0$,
 then the
Boltzmann equation (\ref{(1.1)}) admits a unique global solution
$f^\v(x,t,\xi)$ satisfying
\begin{equation}
\int_{\mb{R}^3}\f{|f^\v(x,t,\xi)-\mb{M}_{[\bar V,\bar
U,\bar\T]}(x,t,\xi)|^2}{\mb{M}_*}d\xi\leq
\tilde{C}_1\delta_0\v^{\f12}+\tilde{C}_2\delta_0e^{-\f{\tilde{C}_3x^2}{\v(1+t)}},
\label{(2.13)}
\end{equation}
with some positive constants $\tilde{C}_i(i=1,2,3)$ independent of
$\v$.

Consequently, we have
\begin{equation}
\sup_{|x|\geq h}\|f^\v(x,t,\xi)-\mb{M}_{[\bar V,\bar
U,\bar\T]}(x,t,\xi)\|_{L_\xi^2(\f{1}{\sqrt{\mb{M}_*}})}\leq C_h
\delta_0 \v^\f14,\qquad \forall h>0,\label{(2.14)}
\end{equation}
where the norm $\|\cdot\|_{L_\xi^2(\f{1}{\sqrt{\mb{M}_*}})}$ is
 $\|\f{\cdot}{\sqrt{\mb{M}_*}}\|_{L_\xi^2(\mb{R}^3)}$.}

\

\noindent{\bf Remark.} Theorem 2.1 shows that, away from the contact
discontinuity located at $x=0$, for any Knudsen number $\v$, there
exists a unique global solution $f^\v(x,t,\xi)$ of the Boltzmann
equation (\ref{(1.1)}) which tends to  $\mb{M}_{[\bar V,\bar
U,\bar\T]}(x,t,\xi)$ as two global Maxwellians with a jump at $x=0$
when $\v\rightarrow 0$. Moreover, a uniform convergence rate
$\v^{\f{1}{4}}$ in the norm
$L_x^{\i}L_\xi^2(\f{1}{\sqrt{\mb{M}_*}})$ holds.

\section{Reformulated system}
\setcounter{equation}{0}

In this section, we will reformulate the system and introduce a
scaling for the independent variable and the perturbation. Firstly,
we define the scaled independent variables
\begin{equation}
y=\v^{-\f12}x,~~\tau=\v^{-\f12}t.\label{(3.1)}
\end{equation}
Correspondingly, set the scaled perturbation as
\begin{equation}
\begin{array}{l}
\di v(x,t)=\bar v(x,t)+\v^{\f12} \p(y,\tau),\\
\di u(x,t)=\bar u(x,t)+\v^{\f12} \psi(y,\tau),\\
\di \t(x,t)=\bar\t(x,t)+\v^{\f12} \z(y,\tau),\\
\di (\t+\f{|u|^2}{2})(x,t)=(\bar \t+\f{|\bar u|^2}{2})(x,t)+\v^{\f12} \omega(y,\tau),\\
\di \mb{G}(x,t,\xi)=\v^{\f12} \bar
\mb{G}(y,\tau,\xi),\\
\di \T_1(x,t,\xi)=\v^{\f12} \bar \T_1(y,\tau,\xi).
\end{array}
\label{(3.2)}
\end{equation}
We remark that the above scaling transformation plays an important
role in the following proof.

Under this scaling,  the hydrodynamic limit problem is now
transferred into a scaled time-asymptotic stability of the viscous
contact wave to the Boltzmann equation. In fact, this scaling is
suitable for the contact wave because of its parabolic structure.
Notice that the hydrodynamic limit proved by this method is globally
in time unlike the case with shock profile proved in \cite{Yu} which
is locally in time. However, we do not know whether there exists
some  appropriate scaling for the shock profile so that this method
can be applied.

With the above scaling, the proof of Theorem 2.1 will be given
by energy method as \cite{Huang-Xin-Yang} for  the scaled perturbation
$(\p,\psi,\z)(y,\tau)$ and $\bar\mb{G}(y,\tau,\xi)$.

From the construction of the contact wave $(\bar v,\bar u,\bar\t)$, the relation between the viscous contact wave $(\bar
v,\bar u,\bar\t)$ to the Boltzmann equation and the inviscid contact
discontinuity $(\bar V,\bar U,\bar\T)$ is given by (\ref{(2.12)}).
Thus, in order to prove Theorem 2.1, it is sufficient to consider
the convergence of the solution $f(y,\tau,\xi)$ of the Boltzmann
equation to the Maxwellian $\mb{M}_{[\bar v,\bar
u,\bar\t]}(y,\tau,\xi)$ defined by the contact wave $(\bar v,\bar
u,\bar\t)$ as the Knudsen number $\v$ tends to zero.

For this, as in \cite{Huang-Xin-Yang}, we introduce the
following anti-derivative of the perturbation:
\begin{equation}
(\Phi,\Psi,\bar
W)(y,\tau)=\int_{-\i}^y(\p,\psi,\omega)(y^\prime,\tau)dy^\prime.
\label{(3.3)}
\end{equation}
Obviously,
$$
(\Phi,\Psi,\bar W)_y(y,\tau)=(\p,\psi,\omega)(y,\tau).
$$
From (\ref{(1.18)}) and (\ref{(2.8)}), we have the following
system for $(\Phi,\Psi,\bar W)$
\begin{equation}
\left\{
\begin{array}{llll}
\di \Phi_{\tau}-\Psi_{1y}=0,\\
\di \Psi_{1\tau}+\v^{-\f12}(p-\bar p)=\f{4}{3}(\f {\mu(\t)}vu_{1y}
-\f{\mu(\bar{\t})}{\bar{v}}\bar{u}_{1y})-\v^{-\f12}R_1-\int\xi_1^2\bar\T_{1}d\xi,\\
\di
\Psi_{i\tau}=(\f{\mu(\t)}{v}u_{iy}-\f{\mu(\bar{\t})}{\bar{v}}\bar{u}_{iy})
-\int\xi_1\xi_i\bar\T_{1}d\xi,~i=2,3,\\
\di
\bar{W}_{\tau}+\v^{-\f12}(pu_1-\bar{p}\bar{u}_1)=(\f{\lambda(\t)}{v}\t_y
 -\f{\lambda(\bar{\t})}{\bar{v}}\bar{\t}_y)+\f43(\f{\mu(\t)}{v}u_1u_{1y}
 -\f{\mu(\bar{\t})}{\bar{v}}\bar{u}_1\bar{u}_{1y})\\
\quad\di+\sum_{i=2}^3\f{\mu(\t)}{v}u_iu_{iy}-\v^{-\f12}R_2-\int\f12\xi_1|\xi|^2\bar\T_1d\xi,
\end{array}
\right. \label{(3.4)}
\end{equation}
where the error terms $R_i~(i=1,2)$ are given in (\ref{(2.9)}) and
(\ref{(2.10)}).

Introduce a new variable
\begin{equation}
W=\bar W-\bar u_1\Psi_1.\label{(3.5)}
\end{equation}
It follows that
\begin{equation}
\z=W_y-Y,~~~{\rm with}~~ Y=\f12\v^{\f12}|\Psi_y|^2-\bar
u_{1y}\Psi_1. \label{(3.6)}
\end{equation}
By using the new variable $W$ and linearizing the system (\ref{(3.4)}),
we have
\begin{equation}
 \left\{
\begin{array}{llll}
\di \Phi_{\tau}-\Psi_{1y}=0,\\
\di \Psi_{1\tau}-\f{p_+}{\bar{v}}\Phi_y+\f{2}{3\bar{v}}W_y=
\f{4}{3}\v^{\f12}\f{\mu(\bar{\t})}{\bar{v}}\Psi_{1yy}-\int\xi_1^2\bar\T_{1}d\xi+Q_1,\\
\di
\Psi_{i\tau}=\v^{\f12}\f{\mu(\bar{\t})}{\bar{v}}\Psi_{iyy}-\int\xi_1\xi_i\bar\T_1d\xi
+Q_i,~i=2,3,\\
\di
 W_{\tau}+p_+\Psi_{1y}=\v^{\f12}\f{\lambda(\bar{\t})}{\bar{v}}W_{yy}
 -\int\f12\xi_1|\xi|^2\bar\T_1d\xi+\bar{u}_1\int\xi_1^2\bar\T_1d\xi+Q_4,
\end{array}
\right. \label{(3.7)}
\end{equation}
where
\begin{equation}
\begin{array}{l}
\di Q_1=\f43(\f {\mu(\t)}v
-\f{\mu(\bar{\t})}{\bar{v}}){u}_{1y}+J_1+\f{2}{3\bar{v}}Y-\v^{-\f12}R_1,\\
\di Q_i=(\f{\mu(\t)}{v}-\f{\mu(\bar{\t})}{\bar{v}}){u}_{iy},~i=2,3,\\
\di Q_4=(\f{\lambda(\t)}{v}
 -\f{\lambda(\bar{\t})}{\bar{v}}){\t}_y
+\f{4\v^\f12}{3}\f{\mu({\t})}{{v}}u_{1y}\Psi_{1y}-\v^{-\f12}R_2-\bar{u}_{1\tau}\Psi_1+\v^{-\f12}\bar{u}_1R_1\\
\di
\qquad+\sum_{i=2}^3\f{\mu(\t)}{v}u_iu_{iy}+J_2-\v^{\f12}\f{\lambda(\bar{\t})}{\bar{v}}Y_y,
\end{array}
\label{(3.8)}
\end{equation}
and
\begin{equation}
\begin{array}{l}
\di J_1=\f{p-p_+}{\bar{v}}\Phi_y
=O(1)(\v^{\f12}\Phi_y^2+\v^{\f12} W_y^2+\v^{\f12} Y^2+|\bar{u}|^4),\\
\di J_2=(p_+-p)\Psi_{1y}=O(1)(\v^{\f12}\Phi_y^2+\v^{\f12}
W_y^2+\v^{\f12}\Psi_{1y}^2+\v^{\f12} Y^2+|\bar{u}|^4).
\end{array}
\label{(3.9)}
\end{equation}
We now derive the equation for the scaled non-fluid component
$\bar\mb{G}(y,\tau,\xi)$. From (\ref{(1.15)}), we have
\begin{equation}
\di
\bar\mb{G}_{\tau}-\f{u_1}{v}\bar\mb{G}_y+\v^{-\f12}\f{1}{v}\mb{P}_1(\xi_1\mb{M}_y)
+\f{1}{v}\mb{P}_1(\xi_1\bar\mb{G}_y)=\v^{-\f12}\mb{L}_\mb{M}\bar\mb{G}
+Q(\bar\mb{G},\bar\mb{G}).\label{(3.10)}
\end{equation}
Thus, we obtain
\begin{equation}
\bar\mb{G}=\f{1}{v}\mb{L}^{-1}_\mb{M}[ \mb{P}_1(\xi_1
\mb{M}_y)]+\bar\T_1,\label{(3.11)}
\end{equation}
and
\begin{equation}
\bar\T_1(y,\tau,\xi)=\v^{\f12}\mb{L}_\mb{M}^{-1}[\bar\mb{G}_{\tau}-\f{u_1}v\bar\mb{G}_y+\f{1}{v}\mb{P}_1(\xi_1\bar\mb{G}_y)
-Q(\bar\mb{G},\bar\mb{G})].\label{(3.12)}
\end{equation}
Let
\begin{equation}
\bar\mb{G}_0(y,\tau,\xi)=\f{3}{2v\t}\mb{L}^{-1}_\mb{M}\{\mb{P}_1[\xi_1(\f{|\xi-u|^2}{2\t}{\bar
\t}_y+\xi\cdot{\bar u}_{y})\mb{M}]\}, \label{(3.13)}
\end{equation}
and
\begin{equation}
\bar\mb{G}_1(y,\tau,\xi)=\bar\mb{G}(y,\tau,\xi)-\bar\mb{G}_0(y,\tau,\xi).
\label{(3.14)}
\end{equation}
Then $\bar\mb{G}_1(y,\tau,\xi)$ satisfies
\begin{equation}
\begin{array}{ll}
\bar\mb{G}_{1\tau}-\v^{-\f12}\mb{L}_\mb{M}\bar\mb{G}_1=&\di-\f{3}{2v\t}\mb{P}_1[\xi_1
(\f{|\xi-u|^2}{2\t}\z_y+\xi\cdot\psi_y)\mb{M}]\\
&\di+\f{u_1}{v}\bar\mb{G}_y-\f{1}{v}\mb{P}_1(\xi_1\bar\mb{G}_y)+Q(\bar\mb{G},\bar\mb{G})-\bar\mb{G}_{0\tau}.
\end{array}
\label{(3.15)}
\end{equation}
Notice that in (\ref{(3.14)}) and (\ref{(3.15)}), $\bar\mb{G}_0$ is
substracted from $\bar\mb{G}$ because $\|\bar{\t}_y\|^2\sim
(1+\v^{\f12}\tau)^{-1/2}$ is not integrable with respect to $\tau$.

Finally, from (\ref{(1.12)}) and the scaling (\ref{(3.1)}), we have
\begin{equation}
\di f_{\tau}-\f{u_1}{v}f_y+\f{\xi_1}{v}f_y=\v^{-\f12}Q(f,f).
\label{(3.16)}
\end{equation}
In the following,
 we will derive the energy estimate on the scaled Boltzmann
equation (\ref{(3.16)}). Indeed, to prove Theorem 2.1, it is sufficient to prove the
following theorem.

\

\noindent{\bf Theorem 3.1.} {\it There exist small positive
constants $\delta_1$, $\v_1$ and a global Maxwellian
$\mb{M}_*=\mb{M}_{[v_*,u_*,\theta_*]}$ such that if the initial data
and the wave strength $\delta$ satisfy
\begin{equation}
E_6(\tau)|_{\tau=0}+\delta \le \delta_1, \label{(3.17)}
\end{equation}
and the Knudsen number $\v$ satisfies $\v\leq \v_1$, then the
problem (\ref{(3.16)}) admits a unique global solution
$f^\v(y,\tau,\xi)$ satisfying
\begin{equation}
\begin{array}{l}
\di
\sup_y\|f^\v(y,\tau,\xi)-\mb{M}_{[\bar{v},\bar{u},\bar{\theta}]}(y,\tau,\xi)
\|_{L^2_{\xi}(\frac{1}{\sqrt{\mb{M}_*}})}\le
C\delta_1\v^{\f12}.\\
\end{array}
\label{(3.18)}
\end{equation}
Here $E_6(\tau)$ will be defined in (\ref{(5.1)}) satisfying
\begin{equation}
\begin{array}{ll}
\di E_6(\tau)\sim&\di
\|(\Phi,\Psi,W)\|^2+\|(\p,\psi,\z)\|^2+\v\|(\p_y,\psi_y,\z_y)\|^2+\int\int\f{\bar\mb{G}_1^2}{\mb{M}_*}d\xi dy\\
&\di +\v\sum_{|\a^\prime|=1} \int\int\f{|\partial^{\a^\prime}
\bar\mb{G}|^2}{\mb{M}_*}d\xi
dy+\v\sum_{|\a|=2}\int\int\f{|\partial^\a f|^2}{\mb{M}_*}d\xi dy.
\end{array}
\label{(3.19)}
\end{equation}}
From now on,  $\partial^\a,\partial^{\a^\prime}$
denote the derivatives with respect to $y$ or $\tau$, and
 $\|\cdot\|^2$ represents $\|\cdot\|_{L_2}^2$ for simplicity of notations.

\

 \noindent{\bf Remark:} In particular,
if we choose the initial value of the Boltzmann equation
(\ref{(3.16)}) as
\begin{equation}
f^\v(y,0,\xi)=\mb{M}_{[\bar v,\bar u,\bar\t]}(y,0,\xi)=\mb{M}_{[\bar
v(y,0),\bar u(y,0),\bar\t(y,0)]}(\xi), \label{(3.20)}
\end{equation}
then
\begin{equation}
E_6(\tau)|_{\tau=0}=O(1)\bigg[\|(\bar \t_y,\bar u_y)\|^2+\v\|(\bar
v_{yy},\bar\t_{yy},\bar u_{yy})\|^2\bigg]\bigg|_{\tau=0}=O(1)\d.
\label{(3.21)}
\end{equation}
In fact, the initial data $f(y,0,\xi)$ can be chosen such that the
initial perturbation $E_6(\tau)|_{\tau=0}$ is suitably small and of
order $O(1)$ with respect to $\v$. This is  the reason why we use
the scaled variables $y,\tau$ in (\ref{(3.1)}), otherwise, the
initial perturbation $E_6(\tau)|_{\tau=0}$ is not uniform with
respect to $\v$.

\section{A priori estimate}
\setcounter{equation}{0}
We will focus on the reformulated system
(\ref{(3.7)}) and (\ref{(3.15)}). Since the local existence of the
solution to (\ref{(3.7)}) and (\ref{(3.15)}) is now standard,
cf. \cite{Ukai-1974} or \cite{Yang-Zhao}, to prove the
global existence,  we
only need to close the following a priori estimate
by the continuity argument
\begin{equation}
\begin{array}{ll}
N(T)=&\di\sup_{0\leq \tau\leq
T}\Bigg\{\|(\Phi,\Psi,W)\|_{L^\i}^2+\|(\p,\psi,\z)\|^2+\v\|(\p_y,\psi_y,\z_y)\|^2+\int\int\f{\bar\mb{G}_1^2}{\mb{M}_*}d\xi
dy\\
&\di+\v\sum_{|\a^\prime|=1} \int\int\f{|\partial^{\a^\prime}
\bar\mb{G}|^2}{\mb{M}_*}d\xi
dy+\v\sum_{|\a|=2}\int\int\f{|\partial^\a f|^2}{\mb{M}_*}d\xi
dy\Bigg\}\leq \gamma^2,
\end{array}
\label{(4.1)}
\end{equation}
where $\gamma$ is a small positive  constant depending on  the
initial data and the strength of the contact wave, and $\mb{M}_*$ is
a global Maxwellian chosen later.

We now briefly explain the a priori assumption
$\|(\Phi,\Psi,W)\|_{L^\i}^2\leq \g^2$ in (\ref{(4.1)}). Roughly
speaking,  based on the observation in \cite{Huang-Xin-Yang} that
the energy estimate involving $\|(\Phi,\Psi,W)\|^2_{L^2}$ may grow
at a rate $(1+\v^{\f12}\tau)^{\frac 12}$, the decay of
$\|(\Phi_x,\Psi_x,W_x)\|^2_{L^2}$  in
 the order of $(1+\v^{\f12}\tau)^{-\frac 12}$
is needed to compensate this growth. This yields a uniform
boundedness of $\|(\Phi,\Psi,W)\|_{L^\i}$, which is essential to
close the a priori estimate.

Note that the a priori assumption (\ref{(4.1)}) also gives
\begin{equation}
\v^{\f12}\|(\p,\psi,\z)\|^2_{L_\i}\leq C\gamma^2, \label{(4.2)}
\end{equation}
\begin{equation}
\v^{\f12}\|\int\f{\bar\mb{G}_1^2}{\mb{M}_*}d\xi\|_{L_\i^y}\leq
C\v^{\f12}\left(\int\int\f{\bar\mb{G}_1^2}{\mb{M}_*}d\xi
dy\right)^{\f{1}{2}}\cdot\left(\int\int\f{|\bar\mb{G}_{1y}|^2}{\mb{M}_*}d\xi
dy\right)^{\f{1}{2}}\leq C(\delta+\gamma)^2, \label{(4.3)}
\end{equation}
and for $|\a^\prime|=1$,
\begin{equation}
\v^{\f32}\|\int\f{|\partial^{\a^\prime}
\bar\mb{G}|^2}{\mb{M}_*}d\xi\|_{L_\i^y}\leq
C\v^{\f32}\left(\int\int\f{|\partial^{\a^\prime}
\bar\mb{G}|^2}{\mb{M}_*}d\xi
dy\right)^{\f{1}{2}}\cdot\left(\int\int\f{|\partial^{\a^\prime}
\bar\mb{G}_y|^2}{\mb{M}_*}d\xi dy\right)^{\f{1}{2}}\leq
C(\delta+\gamma)^2.\label{(4.4)}
\end{equation}
From (\ref{(1.14)}) and (\ref{(2.8)}), we have
\begin{equation}
\left\{
\begin{array}{l}
\di \p_{\tau}-\psi_{1y}=0,\\
\di \psi_{1\tau}+\v^{-\f12}(p-\bar p)_y
              =-\f{4\v^{\f12}}{3}(\f{\mu(\bar\t)}{\bar v}\bar u_{1y})_y-\v^{-\f12}R_{1y}-\int\xi_1^2\bar\mb{G}_yd\xi,\\
\di \psi_{i\tau}=-\v^{-\f12}(\f{\mu(\bar\t)}{\bar v}\bar u_{iy})_y-\int\xi_1\xi_i\bar\mb{G}_yd\xi,~~i=2,3,\\
\di \z_{\tau}+\v^{-\f12}(pu_{1y}-\bar p\bar
u_{1y})=-\v^{-\f12}(\f{\lambda(\bar\t)}{\bar v}\bar \t_y)_y
              -\f{4\v^{\f12}}{3}(\f{\mu(\bar\t)}{\bar v}\bar u_1\bar u_{1y})_y-\v^{-\f12}R_{2y}\\
\di \quad+\v^{-\f12}(\f{|\bar u|^2}{2})_{\tau}-\v^{-\f12}\bar
p_y\bar
u_1-\f{1}{2}\int\xi_1|\xi|^2\bar\mb{G}_yd\xi+\sum_{i=1}^3u_i\int\xi_1\xi_i\bar\mb{G}_yd\xi.
\end{array}
\right. \label{(4.4+)}
\end{equation}
Thus
\begin{equation}
\v\|(\p_{\tau},\psi_{\tau},\z_{\tau})\|^2\leq C(\delta+\gamma)^2.
\label{(4.5)}
\end{equation}
Hence, we have
\begin{equation}
\|(v_{\tau},u_{\tau},\t_{\tau})\|^2\leq
C\v\|(\p_{\tau},\psi_{\tau},\z_{\tau})\|^2+C\|(\bar v_{\tau},\bar
u_{\tau},\bar\t_{\tau})\|^2\leq C(\delta+\gamma)^2. \label{(4.6)}
\end{equation}
In addition,  (\ref{(4.1)}) also implies that
\begin{equation}
\|(v_y,u_y,\t_y)\|^2\leq C\v\|(\p_y,\psi_y,\z_y)\|^2+C\|(\bar
v_y,\bar u_y,\bar\t_y)\|^2\leq C(\delta+\gamma)^2. \label{(4.7)}
\end{equation}
Since
\begin{equation}
\v\|\partial^\a\left(\rho,\rho u,\rho(E+\f{|u|^2}{2})\right)\|^2
\leq C\v\int\int\f{|\partial^\a f|^2}{\mb{M}_*}d\xi dy\le C\gamma^2,
\label{(4.8)}
\end{equation}
(\ref{(4.6)})-(\ref{(4.8)}) give
\begin{equation}
\begin{array}{ll}
\di\v\|\partial^\a(v,u,\t)\|^2&\di\leq
C\v\|\partial^\a\left(\rho,\rho
u,\rho(E+\f{|u|^2}{2})\right)\|^2\\
&\di\quad\quad
+C\v\sum_{|\a^\prime|=1}\int|\partial^{\a^\prime}\left(\rho,\rho
u,\rho(E+\f{|u|^2}{2})\right)|^4dy\\
 &\di\leq C(\delta+\gamma)^2.
\end{array}
\label{(4.9)}
\end{equation}
Thus, for $|\a|=2$, we have
\begin{equation}
\v^2\|\partial^\a(\p,\psi,\z)\|^2\leq
C\v(\|\partial^\a(v,u,\t)\|^2+\|\partial^\a(\bar v,\bar
u,\bar\t)\|^2)\leq C(\delta+\gamma)^2. \label{(4.10)}
\end{equation}
Finally, from the fact that $f=\mb{M}+\v^{\f12}\bar\mb{G}$, we can obtain
for $|\a|=2$,
\begin{equation}
\begin{array}{l}
\di \v^2\int\int\f{|\partial^\a\bar\mb{G}|^2}{\mb{M}_*}d\xi dy\leq
C\v\int\int\f{|\partial^\a f|^2}{\mb{M}_*}d\xi dy+C\v
\int\int\f{|\partial^\a\mb{M}|^2}{\mb{M}_*}d\xi
dy\\
\di\quad \leq C\v\int\int\f{|\partial^\a f|^2}{\mb{M}_*}d\xi
dy+C\v\|\partial^\a(v,u,\t)\|^2+C\v\sum_{|\a^\prime|=1}\int|\partial^{\a^\prime}(v,u,\t)|^4dy\\
\quad \leq C(\delta+\gamma)^2.
\end{array}
\label{(4.11)}
\end{equation}
Before proving the a priori estimate (\ref{(4.1)}), we list some
basic lemmas based on the celebrated H-theorem for later use. The
first lemma is from \cite{Grad}.

\vskip 3mm

 \noindent{\bf Lemma 4.1.} There exists a positive
constant $C$ such that
$$
\int\f{\nu(|\xi|)^{-1}Q(f,g)^2}{\tilde\mb{M}}d\xi\le
C\left\{\int\f{\nu(|\xi|)f^2}{\tilde\mb{M}}d\xi\cdot\int\f{g^2}{\tilde\mb{M}}d\xi+
\int\f{f^2}{\tilde\mb{M}}d\xi\cdot\int\f{\nu(|\xi|)g^2}{\tilde\mb{M}}d\xi\right\},
$$
where $\tilde\mb{M}$ can be any Maxwellian so that the above
integrals are well defined.

\vskip 3mm

Based on Lemma 4.1, the following three lemmas are proved in
\cite{Liu-Yang-Yu-Zhao}. The proofs are straightforward by using
Cauchy inequality.

\vskip 3mm

\noindent{\bf Lemma 4.2.} If $\t/2<\t_*<\t$, then there exist two
positive constants $\sigma=\sigma(v,u,\t;v_*,u_*,\t_*)$ and
$\eta_0=\eta_0(v,u,\t;v_*,u_*,\t_*)$ such that if
$|v-v_*|+|u-u_*|+|\t-\t_*|<\eta_0$, we have for $h(\xi)\in
\mathcal{N}^\bot$,
$$
-\int\f{h\mb{L}_\mb{M}h}{\mb{M}_*}d\xi\geq
\sigma\int\f{\nu(|\xi|)h^2}{\mb{M}_*}d\xi.
$$
\vskip 3mm

 \noindent{\bf Lemma 4.3.} Under the assumptions in Lemma 4.2, we
have  for each $h(\xi)\in \mathcal{N}^\bot$,
$$
\left\{
\begin{array}{l}
\di \int\f{\nu(|\xi|)}{\mb{M}}|\mb{L}_\mb{M}^{-1}h|^2d\xi
\leq \sigma^{-2}\int\f{\nu(|\xi|)^{-1}h^2}{\mb{M}}d\xi,\\
\di \int\f{\nu(|\xi|)}{\mb{M}_*}|\mb{L}_\mb{M}^{-1}h|^2d\xi\le
\sigma^{-2}\int\f{\nu(|\xi|)^{-1}h^2}{\mb{M}_*}d\xi.
\end{array}
\right.
$$

\vskip 3mm

\noindent{\bf Lemma 4.4.} Under the conditions in Lemma 4.2,  for any positive constants $k$ and
$\lambda$, it holds that
$$
|\int\f{g_1\mb{P}_1(|\xi|^kg_2)}{\mb{M}_*}d\xi-\int\f{g_1|\xi|^kg_2}{\mb{M}_*}d\xi|\le
C_{k,\lambda}\int\f{\lambda|g_1|^2+\lambda^{-1}|g_2|^2}{\mb{M}_*}d\xi,
$$
where the constant $C_{k,\lambda}$ depends on $k$ and $\lambda$.

\subsection{Lower order estimate}
Now we will derive the lower
order estimates of $(\Phi,\Psi,W)$. By multiplying $(\ref{(3.7)})_1$ by
$p_+\Phi$, $(\ref{(3.7)})_2$ by $\bar v\Psi_1$, $(\ref{(3.7)})_3$ by
$\Psi_i$, $(\ref{(3.7)})_4$ by $\f{2}{3p_+}W$ respectively and
adding all the resulting equations, we have
\begin{equation}
\begin{array}{l}
\di(\f{p_+}2\Phi^2+\f{\bar{v}}2\Psi_1^2+\f12\sum_{i=2}^3\Psi_i^2+\f{W^2}{3p_+})_{\tau}
+\f{4\v^{\f12}}{3}\mu(\bar{\t})\Psi_{1y}^2+\sum_{i=2}^3\v^{\f12}\f{\mu(\bar{\t})}{\bar{v}}\Psi_{iy}^2
+\f{2\v^{\f12}}{3p_+}\f{\lambda(\bar{\t})}{\bar{v}}W_y^2 \\
=\di\f12\bar{v}_{\tau}\Psi_1^2-\f{4\v^{\f12}}{3}[\mu(\bar{\t})]_y\Psi_{1}\Psi_{1y}
-\sum_{i=2}^3\v^{\f12}(\f{\mu(\bar{\t})}{\bar{v}})_y\Psi_{i}\Psi_{iy}
-\f{2\v^{\f12}}{3p_+}(\f{\lambda(\bar{\t})}{\bar{v}})_yWW_y\\
\di\quad+{\bar v}Q_1\Psi_1
+\sum_{i=2}^3{Q}_i\Psi_i\di+\f{2W}{3p_+}Q_4+N_1 \di+(\cdots)_y,
\end{array}\label{(4.12)}
\end{equation}
where
\begin{equation}
N_1=-\bar{v}\Psi_1\int\xi_1^2\bar\T_1d\xi-\sum_{i=2}^3\Psi_i\int\xi_1\xi_i\bar\T_1d\xi
+\f{2W}{3p_+}(\bar
u_1\int\xi_1^2\bar\T_1d\xi-\int\f12\xi_1|\xi|^2\bar\T_1d\xi).\label{(4.13)}
\end{equation}
From now on,  $(\cdots)_y$ denotes the term
in the conservative form so that it vanishes after integration with
respect to $y$ over $\mb{R}$. Let
\begin{equation}
\begin{array}{l}
\di E_1=\int(\f{p_+}2\Phi^2+\f{\bar{v}}2\Psi_1^2+\f12\sum_{i=2}^3\Psi_i^2+\f{W^2}{3p_+})dy,\\
\di K_1=\int(\f{4\v^{\f12}}{3}\mu(\bar
\t)\Psi_{1y}^2+\sum_{i=2}^3\v^{\f12}\f{\mu(\bar{\t})}{\bar{v}}\Psi_{iy}^2
+\f{2\v^{\f12}}{3p_+}\f{\lambda(\bar{\t})}{\bar{v}}W_y^2)dy.
\end{array}
\label{(4.14)}
\end{equation}
We estimate the right hand side of (\ref{(4.12)}) term by term as follows.
Firstly,
\begin{equation}
\int\f12\bar{v}_{\tau}\Psi_1^2dy\leq
C\delta\v^{\f12}(1+\v^{\f12}\tau)^{-1}E_1, \label{(4.15)}
\end{equation}
\begin{equation}
\int\f{4\v^{\f12}}{3}[\mu(\bar{\t})]_y\Psi_{1}\Psi_{1y}dy\leq \beta
K_1+C_\beta \delta\v^{\f12}(1+\v^{\f12}\tau)^{-1}E_1, \label{(4.16)}
\end{equation}
where $\beta $ is a small positive constant to be
chosen later.

Now we estimate $\int\bar vQ_1\Psi_1dy$ by
\begin{equation}
\begin{array}{ll}
\di \int\bar vQ_1\Psi_1dy
 &\di \leq
\int|\f{4}{3}\bar v(\f{\mu(\t)}{v}-\f{\mu(\bar \t)}{\bar
v})u_{1y}\Psi_1|dy+\int|\bar
vJ_1\Psi_1|dy\\
 &\di+\int|\v^{-\f12}\bar
vR_1\Psi_1|dy+\int|\f{2}{3}Y\Psi_1|dy:=\sum_{i=1}^4 I_i.
\end{array}
\label{(4.17)}
\end{equation}
Note that
\begin{equation}
\begin{array}{l}
I_1 \di \leq C\v^\f12\int|(\Phi_y,\z)\bar
u_{1y}\Psi_1|dy+C\v\int|(\Phi_y,\z)\psi_{1y}\Psi_1|dy\\
    \di\quad\leq
    C\delta\v^{\f12}(1+\v^{\f12}\tau)^{-1}E_1+C(\delta+\gamma)(\v^{\f12}\|\Phi_y\|^2+K_1)+C\gamma\v^{\f32}\|\psi_{1y}\|^2,
\end{array}
\label{(4.18)}
\end{equation}
\begin{equation}
\begin{array}{ll}
I_2 & \di \leq C\int(\v^{\f12}|\Phi_y|^2+\v^{\f12}|W_y|^2+\v^{\f12}
Y^2+|\bar
u|^4)|\Psi_1|dy\\
    & \di\leq
    C\delta\v^{\f12}(1+\v^{\f12}\tau)^{-1}E_1+C\gamma(\v^{\f12}\|\Phi_y\|^2+K_1)+C\delta\v^{\f52}(1+\v^{\f12}\tau)^{-\f{5}{2}},
\end{array}
\label{(4.19)}
\end{equation}
\begin{equation}
I_3\leq C\delta\v^{\f12}(1+\v^{\f12}\tau)^{-1}E_1+C\delta\v^{\f12}
(1+\v^{\f12}\tau)^{-\f{1}{2}}, \label{(4.20)}
\end{equation}
and
\begin{equation}
I_4\leq C\int|(\v^{\f12}|\Psi_y|^2+\bar u_{1y}\Psi_1)\Psi_1| dy\leq
C\delta\v^{\f12}(1+\v^{\f12}\tau)^{-1}E_1+C\gamma K_1.
\label{(4.21)}
\end{equation}
Substituting (\ref{(4.17)})-(\ref{(4.21)}) into (\ref{(4.16)})
yields
\begin{equation}
\int\bar vQ_1\Psi_1dy\leq
C\delta\v^{\f12}(1+\v^{\f12}\tau)^{-1}E_1+C(\delta+\gamma)(\v^{\f12}\|\Phi_y\|^2+K_1)+C\delta
\v^{\f12}(1+\v^{\f12}\tau)^{-\f12}+C\gamma\v^{\f32}\|\psi_{1y}\|^2.
\label{(4.22)}
\end{equation}
Similarly, we can estimate
$$
\int Q_i\Psi_idy~(i=2,3)~~{\rm and}~~\int \f{2W}{3p_+}Q_4dy.
$$
Now we estimate $\int N_1dy$. We only need to estimate
$T_1=:-\int\bar{v}\Psi_1\int\xi_1^2\bar\T_1d\xi dy$ because other
terms in $\int N_1dy$ can be estimated similarly. Let $\mb{M}_*$ be
a global Maxwellian with the state $(v_*,u_*,\t_*)$ satisfying
$\f12\t<\t_*<\t$ and $|v-v_*|+|u-u_*|+|\t-\t_*|\le \eta_0$ so that
Lemma 4.2 holds. By the definition of $\bar\T_1$, cf.
(\ref{(3.12)}), we have
\begin{equation}
\begin{array}{ll}
\di T_1&=\di
-\v^{\f12}\int\bar{v}\Psi_1\int\xi_1^2\mb{L}_\mb{M}^{-1}(\bar\mb{G}_{\tau})d\xi
dy+\v^{\f12}\int\f{u_1\bar{v}\Psi_1}{v}\int\xi_1^2\mb{L}_\mb{M}^{-1}(\bar\mb{G}_y)d\xi
dy\\
&\di\quad
-\v^{\f12}\int\f{\bar{v}\Psi_1}{v}\int\xi_1^2\mb{L}_\mb{M}^{-1}[\mb{P}_1(\xi_1\bar\mb{G}_y)]d\xi
dy+\v^{\f12}\int\bar{v}\Psi_1\int\xi_1^2\mb{L}_\mb{M}^{-1}[Q(\bar\mb{G},\bar\mb{G})]d\xi
dy\\
&\di =: \sum_{i=1}^4T_1^{i}.
\end{array}
\label{(4.23)}
\end{equation}
For the integral $T_1^{1}$, we have
\begin{equation}
\begin{array}{ll}
T_1^{1}&\di
=-\v^{\f12}\int\bar{v}\Psi_1\int\xi_1^2\mb{L}_\mb{M}^{-1}(\bar\mb{G}_{1\tau})d\xi
dy-\v^{\f12}\int\bar{v}\Psi_1\int\xi_1^2\mb{L}_\mb{M}^{-1}(\bar\mb{G}_{0\tau})d\xi
dy\\
&\di =:T_1^{11}+T_1^{12}.
\end{array}
\label{(4.24)}
\end{equation}
Note that the linearized operator $\mb{L}_\mb{M}^{-1}$ satisfies,
for any $h\in \mathcal{N}^\bot$,
\begin{equation}
\begin{array}{l}
(\mb{L}_\mb{M}^{-1}h)_{\tau}=\mb{L}_\mb{M}^{-1}(h_{\tau})-2\mb{L}_\mb{M}^{-1}\{Q(\mb{L}_\mb{M}^{-1}h,M_{\tau})\},\\[2mm]
(\mb{L}_\mb{M}^{-1}h)_y=\mb{L}_\mb{M}^{-1}(h_y)-2\mb{L}_\mb{M}^{-1}\{Q(\mb{L}_\mb{M}^{-1}h,M_y)\}.
\end{array}
\label{(4.25)}
\end{equation}
Then we have
\begin{equation}
\begin{array}{l}
\di
T_1^{11}=-\v^{\f12}\int\bar{v}\Psi_1\int\xi_1^2(\mb{L}_\mb{M}^{-1}\bar\mb{G}_1)_{\tau}d\xi
dy-2\v^{\f12}\int\bar{v}\Psi_1\int\xi_1^2\mb{L}_\mb{M}^{-1}\{Q(\mb{L}_\mb{M}^{-1}\bar\mb{G}_1,M_{\tau})\}d\xi
dy\\
\quad~
\di=-(\v^{\f12}\int\bar{v}\Psi_1\int\xi_1^2\mb{L}_\mb{M}^{-1}(\bar\mb{G}_1)d\xi
dy)_{\tau}+\v^{\f12}\int(\bar{v}\Psi_1)_{\tau}\int\xi_1^2\mb{L}_\mb{M}^{-1}(\bar\mb{G}_1)d\xi
dy\\
\quad\quad~
\di-2\v^{\f12}\int\bar{v}\Psi_1\int\xi_1^2\mb{L}_\mb{M}^{-1}\{Q(\mb{L}_\mb{M}^{-1}\bar\mb{G}_1,\mb{M}_{\tau})\}d\xi
dy.
\end{array}
\label{(4.26)}
\end{equation}
The H\"{o}lder inequality and Lemma 4.3 yield
\begin{equation}
|\int\xi_1^2\mb{L}_\mb{M}^{-1}(\bar\mb{G}_1)d\xi|^2\leq
C\int\f{\nu^{-1}(|\xi|)}{\mb{M}_*}|\bar\mb{G}_1|^2d\xi.
\label{(4.27)}
\end{equation}
Moreover, from Lemmas 4.1-4.3, we have
\begin{equation}
\begin{array}{l}
\di
\int\xi_1^2\mb{L}_\mb{M}^{-1}\{Q(\mb{L}_\mb{M}^{-1}\bar\mb{G}_1,\mb{M}_{\tau})\}d\xi\leq
C\left(\int\f{\nu(|\xi|)}{\mb{M}_*}|\mb{L}_\mb{M}^{-1}\{Q(\mb{L}_\mb{M}^{-1}\bar\mb{G}_1,\mb{M}_{\tau})\}|^2d\xi\right)^{\f{1}{2}}
\\ \quad \di
\leq
C\left(\int\f{\nu(|\xi|)}{\mb{M}_*}|\mb{L}_\mb{M}^{-1}\bar\mb{G}_1|^2d\xi\right)^{\f{1}{2}}\cdot\left(\int\f{\nu(|\xi|)}{\mb{M}_*}|\mb{M}_{\tau}|^2d\xi\right)^{\f{1}{2}}
\\
\quad \di\leq
C|(v_{\tau},u_{\tau},\t_{\tau})|\left(\int\f{\nu^{-1}(|\xi|)}{\mb{M}_*}|\bar\mb{G}_1|^2d\xi\right)^{\f{1}{2}}.
\end{array}
\label{(4.28)}
\end{equation}
Combining (\ref{(4.26)})-(\ref{(4.28)}) gives
\begin{equation}
\begin{array}{ll}
\di T_1^{11}\leq &\di
-(\v^{\f12}\int\bar{v}\Psi_1\int\xi_1^2\mb{L}_\mb{M}^{-1}\bar\mb{G}_1d\xi
dy)_{\tau}+C\delta\v^{\f12}(1+\v^{\f12}\tau)^{-1}E_1\\
&\di+C\beta
\v^{\f12}\|\Psi_{1\tau}\|^2+C\v^{\f12}\int\int\f{\nu(|\xi|)}{\mb{M}_*}|\bar\mb{G}_1|^2d\xi
dy+C\gamma\v^{\f32}\|(\p_{\tau},\psi_{\tau},\z_{\tau})\|^2.
\end{array}
\label{(4.29)}
\end{equation}
On the other hand, by (\ref{(3.13)}), we have
\begin{equation}
\begin{array}{l}
\di
T_1^{12}=-\v^{\f12}\int\bar{v}\Psi_1\int\xi_1^2\mb{L}_\mb{M}^{-1}(\bar\mb{G}_{0\tau})d\xi
dy\\
\quad\di \leq C\v^\f12\int
|\Psi_1|(|(\bar{\t}_{y\tau},\bar{u}_{y\tau})|
+|(\bar{\t}_y,\bar{u}_y)||(v_{\tau},u_{\tau},\t_{\tau})|)dy\\
\quad \di\leq
C\delta\v^{\f12}(1+\v^{\f12}\tau)^{-1}E_1+C\delta\v(1+\v^{\f12}\tau)^{-\f32}+C\delta\v^{\f32}\|(\p_{\tau},\psi_{\tau},
\z_{\tau})\|^2,
\end{array}
\label{(4.30)}
\end{equation}
which, together with (\ref{(4.29)}), imply
\begin{equation}
\begin{array}{l}
\di T_1^{1}\leq
-(\v^{\f12}\int\bar{v}\Psi_1\int\xi_1^2\mb{L}_\mb{M}^{-1}\bar\mb{G}_1d\xi
dy)_{\tau}+C\delta\v^{\f12}(1+\v^{\f12}\tau)^{-1}E_1+C\delta\v(1+\v^{\f12}\tau)^{-\f32}\\
\quad \di+C\beta
\v^{\f12}\|\Psi_{1\tau}\|^2+C\v^{\f12}\int\int\f{\nu(|\xi|)}{\mb{M}_*}|\bar\mb{G}_1|^2d\xi
dy+C(\delta+\gamma)\v^{\f32}\|(\p_{\tau},\psi_{\tau},\z_{\tau})\|^2.
\end{array}
\label{(4.31)}
\end{equation}

The estimation on $T_1^{i}$, $i=2,4$ is relatively easy by using the
Cauchy inequality and Lemmas 4.1-4.3. In fact,
direct computation gives
\begin{equation}
\begin{array}{l}
\di T_1^{2}\leq
C\v^{\f12}\int\int\f{\nu(|\xi|)}{\mb{M}_*}|\bar\mb{G}_y|^2d\xi
dy+C\v^{\f12}\int
\Psi_1^2u_1^2dy\\
\quad\quad\di\leq C\delta\v^2(1+\v^{\f12}\tau)^{-2}E_1+C\gamma\v
K_1+C\v^{\f12}\int\int\f{\nu(|\xi|)}{\mb{M}_*}|\bar\mb{G}_y|^2d\xi
dy.
\end{array}
\label{(4.32)}
\end{equation}
On the other hand,
\begin{equation}
\begin{array}{l}
\di T_1^{4}\leq
C\gamma\v^{\f12}\int(\int\f{\nu(|\xi|)}{\mb{M}_*}|\mb{L}_\mb{M}^{-1}\{Q(\bar\mb{G},\bar\mb{G})\}|^2d\xi)^{\f12}dy
\\
\quad\di \leq
C\gamma\v^{\f12}\int\int\f{\nu(|\xi|)}{\mb{M}_*}|\bar\mb{G}|^2d\xi
dy\\
\quad\di \leq
C\gamma\v^{\f12}\int\int\f{\nu(|\xi|)}{\mb{M}_*}|\bar\mb{G}_1|^2d\xi
dy+C\delta\v^{\f12}(1+\v^{\f12}\tau)^{-\f12}.
\end{array}
\label{(4.33)}
\end{equation}

The estimation on $T_1^{3}$ is similar to the one for $T_1^{1}$.
Firstly, notice that
\begin{equation}
\mb{P}_1(\xi_1\bar\mb{G}_y)=[\mb{P}_1(\xi_1\bar\mb{G})]_y+\sum_{j=0}^4<\xi_1\bar\mb{G},\chi_j>\mb{P}_1(\chi_{jy}).
\label{(4.34)}
\end{equation}
 Then, it follows
from (3.46), (3.55) and Lemmas 3.1-3.4 that
\begin{equation}
\begin{array}{ll}
\di
T_1^{3}=&\di\v^{\f12}\int(\f{\bar{v}\Psi_1}{v})_y\int\xi_1^2\mb{L}_\mb{M}^{-1}[\mb{P}_1(\xi_1\bar\mb{G})]d\xi
dy\\
\quad
&\di-\v^{\f12}\int\f{\bar{v}\Psi_1}{v}\int\xi_1^2\mb{L}_\mb{M}^{-1}
[\sum_{j=0}^4<\xi_1\bar\mb{G},\chi_j>\mb{P}_1(\chi_{jy})]d\xi dy
\\
\quad
&\di-2\v^{\f12}\int\f{\bar{v}\Psi_1}{v}\int\xi_1^2\mb{L}_\mb{M}^{-1}
\{Q(\mb{L}_\mb{M}^{-1}[\mb{P}_1(\xi_1\bar\mb{G})],\mb{M}_y)\}d\xi
dy\\
\quad \leq &\di
C\delta\v^{\f12}(1+\v^{\f12}\tau)^{-1}E_1+C(\gamma+\beta
)K_1+C\delta\v^{\f12}(1+\v^{\f12}\tau)^{-\f12}+C\gamma\v^{\f32}\|(\p_y,\psi_y,\z_y)\|^2\\
\quad &\di
+C\v^{\f12}\int\int\f{\nu(|\xi|)}{\mb{M}_*}|\bar\mb{G}_1|^2d\xi dy.
\end{array}
\label{(4.35)}
\end{equation}
By (\ref{(4.24)}), (\ref{(4.26)})-(\ref{(4.33)}) and (\ref{(4.35)}),
we have
\begin{equation}
\begin{array}{l}
\di T_1\leq
-(\v^{\f12}\int\bar{v}\Psi_1\int\xi_1^2\mb{L}_\mb{M}^{-1}(\bar\mb{G}_1)d\xi
dy)_{\tau}+C\delta\v^{\f12}(1+\v^{\f12}\tau)^{-1}E_1+C\delta\v^{\f12}(1+\v^{\f12}\tau)^{-\f12}\\
\quad\di +C\beta \v^{\f12}\|\Psi_{1\tau}\|^2+C(\gamma+\beta
)(K_1+\v^{\f12}\|\Phi_y\|^2)\di+C\v^{\f12}\int\int\f{\nu(|\xi|)}{\mb{M}_*}|\bar\mb{G}_1|^2d\xi
dy\\
\quad\di+C\v^{\f12}\int\int\f{\nu(|\xi|)}{\mb{M}_*}|\bar\mb{G}_y|^2d\xi
dy+C(\delta
+\gamma)\v^{\f32}\sum_{|\a^\prime|=1}\|\partial^{\a^\prime}(\p,\psi,\z)\|^2.
\end{array}
\label{(4.36)}
\end{equation}

The estimates on the other terms of $\int N_1dy$ are similar and we
omit the details for brevity.
Therefore, collecting the above inequalities gives
\begin{equation}
\begin{array}{l}
\di
E_{1\tau}+(\v^{\f12}\int\int\hat{A}(\xi,\Phi,\Psi,W)\mb{L}_\mb{M}^{-1}(\bar\mb{G}_1)d\xi
dy)_{\tau}+\f{1}{2} K_1\leq
C_1\delta\v^{\f12}(1+\v^{\f12}\tau)^{-1}E_1\\
\di+C_1 \beta
\v^{\f12}\|(\Psi_{\tau},W_{\tau})\|^2+C_1(\delta+\gamma)\v^{\f12}\|\Phi_y\|^2+C_1\v^{\f12}\int\int\f{\nu(|\xi|)}{\mb{M}_*}|\bar\mb{G}_1|^2d\xi
dy\\
\di+C_1\v^{\f12}\int\int\f{\nu(|\xi|)}{\mb{M}_*}|\bar\mb{G}_y|^2d\xi
dy+C_1(\delta+\gamma)\v^{\f32}\sum_{|\a^\prime|=1}\|\partial^{\a^\prime}(\p,\psi,\z)\|^2+C_1\delta\v^{\f12}(1+\v^{\f12}\tau)^{-\f12}
,
\end{array}
\label{(4.37)}
\end{equation}
where we have used the smallness of $\delta$, $\beta$ and $\gamma$.
Here $\hat{A}(\xi,\Phi,\Psi,W)$ is a linear function of
$(\Phi,\Psi,W)$ which is  a polynomial of $\xi$.

\

Note that the dissipation term
$K_1$ does not contain the term $\v^{\f12}\|\Phi_y\|^2$.
To complete the lower order inequality, we have to estimate
$\Phi_y$. From $(3.8)_2$, we have
\begin{equation}
\f{4\v^{\f12}}{3}\f{\mu(\bar{\t})}{\bar{v}}\Phi_{y\tau}-\Psi_{1\tau}+\f{p_+}{\bar{v}}\Phi_y=\f2{3\bar{v}}W_y-Q_1+\int\xi_1^2\bar\T_1d\xi.
\label{(4.38)}
\end{equation}
Multiplying (\ref{(4.38)}) by $\v^{\f12}\Phi_y$ yields
\begin{equation}
\begin{array}{l}
\quad\di(\v\f{2\mu(\bar{\t})}{3\bar{v}}\Phi_y^2-\v^{\f12}\Phi_y\Psi_1)_{\tau}+\v^{\f12}
\f{p_+}{\bar{v}}\Phi_y^2\\
\di
=\v(\f{2\mu(\bar\t)}{3\bar{v}})_{\tau}\Phi_y^2+\v^{\f12}\Psi_{1y}^2+\v^{\f12}(\f{2}{3\bar{v}}W_y-Q_1+\int\xi_1^2\bar\T_1d\xi)\Phi_y+(\cdots)_y,
\end{array}
\label{(4.39)}
\end{equation}
where we have used
$$
\Phi_y\Psi_{1\tau}=(\Phi_y\Psi_1)_{\tau}
-(\Phi_{\tau}\Psi_1)_y+\Psi_{1y}^2.
$$
Integrating (\ref{(4.39)}) with respect to $y$ gives
\begin{equation}
\begin{array}{l}
\di
(\int\v\f{2\mu(\bar{\t})}{3\bar{v}}\Phi_y^2-\v^{\f12}\Phi_y\Psi_1dy)_{\tau}+\int\v^{\f12}\f{p_+}{2\bar{v}}\Phi_y^2dy\\
\quad\di\leq CK_1+C\v^{\f12}\int
Q_1^2dy+C\v^{\f12}\int|\int\xi_1^2\bar\T_1d\xi|^2dy.
\end{array}
\label{(4.40)}
\end{equation}
By (\ref{(3.8)}) and the Cauchy inequality, one has
\begin{equation}
\v^{\f12}\int Q_1^2dy\leq C\delta\v^{\f12}(1+\v^{\f12}\tau)^{-1}E_1+
C\gamma\v(K_1+\v^{\f12}\|\Phi_y\|^2)+C\delta\v^{\f32}(1+\v^{\f12}\tau)^{-\f32}+
C(\delta+\gamma)\v^2\|\psi_{1y}\|^2. \label{(4.41)}
\end{equation}
On the other hand, Lemmas 4.1-4.3 imply
\begin{equation}
\begin{array}{ll}
\di\v^{\f12}\int|\int\xi_1^2\bar\T_1d\xi|^2dy &\di\leq
C\v^{\f32}\sum_{|\a^\prime|=1}\int\int\f{\nu(|\xi|)}{\mb{M}_*}|\partial^{\a^\prime}\bar\mb{G}|^2d\xi
dy\\
&\quad\di+C\gamma\v\int\int\f{\nu(|\xi|)}{\mb{M}_*}|\bar\mb{G}_1|^2d\xi
dy+C\delta\v^{\f32}(1+\v^{\f12}\tau)^{-\f32}.
\end{array}
\label{(4.42)}
\end{equation}
Thus combining (\ref{(4.41)})-(\ref{(4.42)}) yields
\begin{equation}
\begin{array}{l}
\di
(\int\v\f{2\mu(\bar{\t})}{3\bar{v}}\Phi_y^2-\v^{\f12}\Phi_y\Psi_1dy)_{\tau}+\int\v^{\f12}\f{p_+}{4\bar{v}}\Phi_y^2dy\\
\quad\di\leq
C_2\delta\v^{\f12}(1+\v^{\f12}\tau)^{-1}E_1+C_2K_1+C_2\delta\v^{\f32}(1+\v^{\f12}\tau)^{-\f32}+
C(\delta+\gamma)\v^2\|\psi_{1y}\|^2\\
\quad\di
+C_2\v^{\f32}\sum_{|\a^\prime|=1}\int\int\f{\nu(|\xi|)}{\mb{M}_*}|\partial^{\a^\prime}\bar\mb{G}|^2d\xi
dy+C_2\gamma\v\int\int\f{\nu(|\xi|)}{\mb{M}_*}|\bar\mb{G}_1|^2d\xi
dy.
\end{array}
\label{(4.43)}
\end{equation}

The microscopic component $\bar\mb{G}_1$ can be estimated through
the equation (\ref{(3.15)}). Multiplying (\ref{(3.15)}) by
$\f{\bar\mb{G}_1}{\mb{M}_*}$ gives
\begin{equation}
\begin{array}{ll}
(\f{\bar\mb{G}_1^2}{2\mb{M}_*})_{\tau}-\v^{-\f12}\f{\bar\mb{G}_1}{\mb{M}_*}\mb{L}_\mb{M}\bar\mb{G}_1=&\di\bigg\{-\f1{Rv\t}\mb{P}_1[\xi_1
(\f{|\xi-u|^2}{2\t}\z_y+\xi\cdot\psi_y)\mb{M}]\\
&\di+\f{u_1}{v}\bar\mb{G}_y-\f1v\mb{P}_1(\xi_1\bar\mb{G}_y)+Q(\bar\mb{G},\bar\mb{G})-\bar\mb{G}_{0\tau}\bigg\}
\f{\bar\mb{G}_1}{\mb{M}_*}.
\end{array}
\label{(4.44)}
\end{equation}
Integrating (\ref{(4.44)}) with respect to $\xi$ and $y$ and using
the Cauchy inequality and Lemma 4.1-4.4 yield that
\begin{equation}
\begin{array}{l}
\di(\int\int\f{\bar\mb{G}_1^2}{2\mb{M}_*}d\xi
dy)_{\tau}+\f{\sigma}{2}\v^{-\f12}\int\int\f{\nu(|\xi|)}{\mb{M}_*}|\bar\mb{G}_1|^2d\xi dy\\
\leq\di
C_3\delta\v^{\f12}(1+\v^{\f12}\tau)^{-\f12}+C_3\v^{\f12}\sum_{|\a^\prime|=1}\|\partial^{\a^\prime}(\p,\psi,\z)\|^2
+C_3\v^{\f12}\int\int\f{\nu(|\xi|)}{\mb{M}_*}|\bar\mb{G}_y|^2d\xi
dy.
\end{array}
\label{(4.45)}
\end{equation}

On the other hand, from the fluid-type system  (\ref{(3.7)}), we
can get an estimate for $\v^{\f12}\|(\Psi_{\tau},W_{\tau})\|^2$ as
follows.
\begin{equation}
\begin{array}{ll}
\di \v^{\f12}\|(\Psi_{\tau},W_{\tau})\|^2\leq&\di
C_4\v^{\f12}(1+\v^{\f12}\tau)^{-1}E_1+C_4K_1+C_4\v^{\f12}\|\Phi_y\|^2+C_4\v^{\f32}\|(\psi_y,\z_y)\|^2\\
&\di+C_4\delta\v^{\f32}(1+\v^{\f12}\tau)^{-\f32}+C_4(\delta+\gamma)\v^{\f12}\int\int\f{\nu(|\xi|)
}{\mb{M}_*}|\bar\mb{G}_1|^2d\xi
dy\\
&\di+C_4\v^{\f32}\sum_{|\a^\prime|=1}\int\int\f{\nu(|\xi|)}{\mb{M}_*}|\partial^{\a^\prime}\bar\mb{G}|^2d\xi
dy .
\end{array}
\label{(4.46)}
\end{equation}
We can now complete the lower order estimate. Since
$\hat{A}(\xi,\Phi,\Psi,W)$ is a linear function of the vector
$(\Phi,\Psi,W)$ which is a polynomial of $\xi$, we get
$$
|\v^{\f12}\int\int\hat{A}\mb{L}_\mb{M}^{-1}\bar\mb{G}_1d\xi dy|\leq
\f14E_1+C\v\int\int\f{\bar\mb{G}_1^2}{\mb{M}_*}d\xi dy.
$$
We choose large constants $\bar{C}_1>1$, $\bar{C}_2>1$,
$\bar{C}_3>1$ and small constant $\beta $ such that
\begin{equation}
\begin{array}{lll}
\di E_2&=&\di
\bar{C}_1E_1+\bar{C}_1\v^{\f12}\int\int\hat{A}\mb{L}_\mb{M}^{-1}\bar\mb{G}_1d\xi
dy+\bar{C}_2\int\v\f{2\mu(\bar{\t})}{3\bar{v}}\Phi_y^2-\v^{\f12}\Phi_y\Psi_1dy\\
&&\di+ \bar{C}_3\int\int\f{\bar\mb{G}_1^2}{2\mb{M}_*}d\xi dy \\
&\geq&\di
\f12\bar{C}_1E_1+\bar{C}_2\int\v\f{\mu(\bar{\t})}{3\bar{v}}\Phi_y^2dy
+\f{\bar{C}_3}{4}\int\int\f{\bar\mb{G}_1^2}{\mb{M}_*}d\xi dy,
\end{array}
\label{(4.47)}
\end{equation}
and
\begin{equation}
\begin{array}{l}
\di
 (\f{\bar{C}_1}{2}-C_2\bar{C}_2
 -\bar{C}_1C_1\beta C_4)K_1+\int\v^{\f12}(\bar{C}_2\f{p_+}{4\bar{v}}
 -\bar{C}_1C_1\beta (1+C_4))\Phi_y^2dy\\
\geq\di  \f{\bar{C}_1}{4}K_1+
\bar{C}_2\int\v^{\f12}\f{p_+}{8\bar{v}}\Phi_y^2dy.
\end{array}
\label{(4.48)}
\end{equation}
Hence, multiplying (\ref{(4.37)}) by $\bar{C}_1$, (\ref{(4.43)}) by
$\bar{C}_2$, (\ref{(4.45)}) by $\bar{C}_3$, (\ref{(4.46)}) by
$C_1(\delta+\gamma)\bar{C}_1$ and adding all these
inequalities imply
\begin{equation}
\begin{array}{ll}
\di E_{2\tau}+K_2\leq&\di
C_5\delta\v^{\f12}(1+\v^{\f12}\tau)^{-1}E_2+C_5\v^{\f12}\sum_{|\a^\prime|=1}\int\int\f{\nu(|\xi|)}{\mb{M}_*}|\partial^{\a^\prime}\bar\mb{G}|^2d\xi dy\\
&\di+C_5\v^{\f12}\sum_{|\a^\prime|=1}\|\partial^{\a^\prime}(\p,\psi,\z)\|^2
+C_5\delta\v^{\f12}(1+\v^{\f12}\tau)^{-\f12},
\end{array}\label{(4.49)}
\end{equation}
where
\begin{equation}
K_2=\f{\bar{C}_1}{4}K_1+\bar{C}_2
\int\v^{\f12}\f{p_+}{8\bar{v}}\Phi_y^2dy+\v^{\f12}\|(\Psi_\tau,W_\tau)\|^2+\f{\sigma}{4}\bar{C}_3\v^{-\f12}\int\int\f{\nu(|\xi|)}{\mb{M}_*}|\bar\mb{G}_1|^2d\xi
dy.\label{(4.50)}
\end{equation}

\subsection{Higher order estimate}
In this subsection, we shall
estimate the derivatives of $(\Phi,\Psi,W)$. Applying $\partial_y$
to the system (\ref{(3.4)}) gives
\begin{equation}
 \left\{
\begin{array}{l}
\di \p_{\tau}-\psi_{1y}=0,\\
\di \psi_{1\tau}+\v^{-\f12}(p-\bar{p})_y=
\f{4}{3}(\f{\mu(\t)}{v}u_{1y}-\f{\mu(\bar\t)}{\bar{v}}\bar{u}_{1y})_y-\v^{-\f12} R_{1y}-\int\xi_1^2\bar\T_{1y}d\xi,\\
\di \psi_{i\tau}=
(\f{\mu(\t)}{v}u_{iy}-\f{\mu(\bar\t)}{\bar{v}}\bar{u}_{iy})_y-\int\xi_1\xi_i\bar\T_{1y}d\xi,~~i=2,3,\\
\di \z_{\tau}+\v^{-\f12}(pu_{1y}-\bar{p}\bar{u}_{1y})
=(\f{\lambda(\t)}{v}\t_{y}
-\f{\lambda(\bar\t)}{\bar{v}}\bar{\t}_{y})_y+Q_5\\
\hspace{3cm}\di+\sum_{i=1}^3u_i\int\xi_1\xi_i\bar\T_{1y}d\xi-\f12\int\xi_1|\xi|^2\bar\T_{1y}d\xi,
\end{array}
\right. \label{(4.51)}
\end{equation}
where
\begin{equation}
\begin{array}{ll}
Q_5&=\di\f{4}{3}(\f{\mu(\t)}{v}u_{1y}^2-\f{\mu(\bar\t)}{\bar v}\bar
u_{1y}^2)+\sum_{i=2}^3\f{\mu(\t)}{v}u_{iy}^2-\v^{-\f12}R_{2y}-\v^{-\f12}R
_{1y}\bar u_1\\
&\di =O(1)\left[\v |\psi_y|^2+\v^{\f12}|\bar
u_{1y}|^2|(\p,\z)|+\v^\f12|\psi_{1y}||\bar
u_{1y}|+\v^{-\f12}|R_{2y}|+\v^{-\f12}|R_{1y}\bar u_1|\right].
\end{array}\label{(4.52)}
\end{equation}
Multiplying $(\ref{(4.51)})_2$ by $\psi_1$, $(\ref{(4.51)})_3$ by
$\psi_i~(i=2,3)$ respectively and adding them together yield
$$
\begin{array}{ll}
\di(\sum_{i=1}^3\f12\psi_i^2)_{\tau}-\v^{-\f12}(p-\bar{p})\psi_{1y}+\f43(\f{\mu(\t)}{v}u_{1y}
-\f{\mu(\bar\t)}{\bar{v}}\bar{u}_{1y})\psi_{1y}\\
\quad\di+\sum_{i=2}^3(\f{\mu(\t)}{v}u_{iy}
-\f{\mu(\bar\t)}{\bar{v}}\bar{u}_{iy})\psi_{iy}=-\v^{-\f12}R_{1y}\psi_1
-\sum_{i=1}^3\psi_i\int\xi_1\xi_i\bar\T_{1y}d\xi+(\cdots)_y.
\end{array}
$$
Since $p-\bar{p}=\v^{\f12}\f{R\z}{v}+R\bar{\t}(\f1v-\f1{\bar v}),$
we obtain
\begin{equation}
\begin{array}{l}
\di(\sum_{i=1}^3\f12\psi_i^2)_{\tau}-\v^{-\f12}R\bar{\t}(\f1v-\f1{\bar
v})\p_{\tau}-\f{R\z}{v}\psi_{1y}+\f43\v^{\f12}\f{\mu(\t)}{v}\psi_{1y}^2
+\sum_{i=2}^3\v^{\f12}\f{\mu(\t)}{v}\psi_{iy}^2\\
\quad\di=-\f43(\f{\mu(\t)}{v}
-\f{\mu(\bar\t)}{\bar{v}})\bar{u}_{1y}\psi_{1y}-\v^{-\f12}R_{1y}\psi_1
-\sum_{i=1}^3\psi_i\int\xi_1\xi_i\bar\T_{1y}d\xi+(\cdots)_y.
\end{array}
\label{(4.53)}
\end{equation}
Set
$$
\hat{\Phi}(s)=s-1-\ln s.
$$
Then
\begin{equation}
\begin{array}{l}
\di \{R\bar{\t}\hat{\Phi}(\f{v}{\bar{v}})\}_{\tau}=-\v^{\f12}
R\bar{\t}(\f{1}{v}-\f{1}{\bar{v}})\p_{\tau}-\bar{p}\hat{\Psi}(\f{v}{\bar{v}})\bar{v}_{\tau}
+\bar{v}\bar{p}_{\tau}\hat{\Phi}(\f{v}{\bar{v}}),
\end{array}
\label{(4.54)}
\end{equation}
where
$$
\hat{\Psi}(s)=s^{-1}-1+\ln{s}.
$$
It is easy to check that
$\hat{\Phi}(1)=\hat{\Phi}'(1)=\hat{\Psi}(1)=\hat{\Psi}'(1)=0$ and
$\hat{\Phi}(s)$ is strictly convex around $s=1$. Substituting
(\ref{(4.54)}) into (\ref{(4.53)}) yields
\begin{equation}
\begin{array}{l}
\di(\sum_{i=1}^3\f12\psi_i^2+\v^{-1}R\bar{\t}\hat{\Phi}(\f{v}{\bar{v}}))_{\tau}
-\f{R}{v}\z\psi_{1y}+\f43\v^{\f12}\f{\mu(\t)}{v}\psi_{1y}^2
+\sum_{i=2}^3\v^{\f12}\f{\mu(\t)}{v}\psi_{iy}^2=-\v^{-1}\bar{p}\hat{\Psi}(\f{v}{\bar{v}})\bar{v}_{\tau}\\
\quad\di
+\v^{-1}\bar{v}\bar{p}_{\tau}\hat{\Phi}(\f{v}{\bar{v}})-\f43(\f{\mu(\t)}{v}
-\f{\mu(\bar\t)}{\bar{v}})\bar{u}_{1y}\psi_{1y}-\v^{-\f12}R_{1y}\psi_1
-\sum_{i=1}^3\psi_i\int\xi_1\xi_i\bar\T_{1y}d\xi+(\cdots)_y.
\end{array}
\label{(4.55)}
\end{equation}
Note that
\begin{equation}
\{\bar{\t}
\hat{\Phi}(\f{\t}{\bar{\t}})\}_{\tau}=\v^{\f12}(1-\f{\bar{\t}}{\t})\z_{\tau}-
\hat{\Psi}(\f{\t}{\bar{\t}})\bar{\t}_{\tau},\label{(4.56)}
\end{equation}
and
\begin{equation}
\begin{array}{ll}
&\di \v^{\f12}(1-\f{\bar{\t}}{\t})\z_{\tau}\\
&\di
  =\v\f{\z}{\t}\bigg\{-\v^{-\f12}(pu_{1y}-\bar{p}\bar{u}_{1y})
  +(\f{\lambda(\t)}{v}\t_{y}
  -\f{\lambda(\bar\t)}{\bar{v}}\bar{\t}_{y})_y+Q_5\\
&\di \qquad\qquad +\sum_{i=1}^3u_i\int\xi_1\xi_i\bar\T_{1y}d\xi-\f12\int\xi_1|\xi|^2\bar\T_{1y}d\xi\bigg\}\\
\di
  &\di =-\v\f{R\z}{v}\psi_{1y}+\v\bigg\{\v^{-\f12}\f{\z}{\t}(\bar{p}-p)\bar{u}_{1y}
  -\v^{\f12}\f{\lambda(\t)}{v\t}\z_y^2-\f{\z_y}{\t}(\f{\lambda(\t)}{v}-\f{\lambda(\bar\t)}{\bar{v}})\bar{\t}_y
  +\f{\z}{\t}Q_5\\
&\di\quad \di+\f{\z\t_y}{\t^2}(\f{\lambda(\t)
  \t_y}{v}-\f{\lambda(\bar\t)\bar{\t}_y}{\bar{v}})+\f{\z}{\t}(\sum_{i=1}^3u_i\int\xi_1\xi_i\bar\T_{1y}d\xi-\f12\int\xi_1|\xi|^2\bar\T_{1y}d\xi)+(\cdots)_y\bigg\}.
\end{array}
\label{(4.57)}
\end{equation}
Substituting (\ref{(4.56)}) and (\ref{(4.57)}) into (\ref{(4.55)})
gives
\begin{equation}
\begin{array}{l}
\di\left(\sum_{i=1}^3\f12\psi_i^2+\v^{-1}R\bar{\t}\hat{\Phi}(\f{v}{\bar{v}})+\v^{-1}\bar{\t}
\hat{\Phi}(\f{\t}{\bar{\t}})\right)_{\tau}
+\f{4\v^{\f12}}{3}\f{\mu(\t)}{v}\psi_{1y}^2
+\sum_{i=2}^3\v^{\f12}\f{\mu(\t)}{v}\psi_{iy}^2+\v^{\f12}\f{\lambda(\t)}{v\t}\z_y^2\\
\di=-\v^{-1}\bar{p}\hat{\Psi}(\f{v}{\bar{v}})\bar{v}_{\tau}+\v^{-1}\bar{v}\bar{p}_{\tau}\hat{\Phi}(\f{v}{\bar{v}})-\v^{-1}\hat{\Psi}(\f{\t}{\bar{\t}})\bar{\t}_{\tau}
-\f{4}{3}(\f{\mu(\t)}{v}-\f{\mu(\bar{\t})}{\bar{v}})
\bar{u}_{1y}\psi_{1y}
+\v^{-\f12}R_{1y}\psi_1\\[0.3cm]
\di+\v^{-\f12}\f{\z}{\t}(\bar{p}-p)\bar{u}_{1y}
-\f{\z_y}{\t}(\f{\lambda(\t)}{v}
-\f{\lambda(\bar{\t})}{\bar{v}})\bar{\t}_y+\f{\z\t_y}{\t^2}(\f{\lambda(\t)\t_y}{v}-\f{\lambda(\bar{\t})\bar{\t}_y}{\bar{v}})
+\f{\z}{\t}Q_5+N_2+(\cdots)_y,
\end{array}
\label{(4.58)}
\end{equation}
where
\begin{equation}
N_2=-\sum_{i=1}^3\psi_i\int\xi_1\xi_i\bar\T_{1y}d\xi+
\f{\z}{\t}(\sum_{i=1}^3u_i\int\xi_1\xi_i\bar\T_{1y}d\xi-\f12\int\xi_1|\xi|^2
\bar\T_{1y}d\xi). \label{(4.59)}
\end{equation}
Let
\begin{equation}
\begin{array}{l}
\di
E_3=\int(\f12\sum_{i=1}^3\psi_i^2+\v^{-1}R\bar{\t}\hat{\Phi}(\f{v}{\bar{v}})+\v^{-1}\bar{\t}
\hat{\Phi}(\f{\t}{\bar{\t}}))dy,\\
\di K_3=\int(\f{4\v^{\f12}}{3}\f{\mu(\t)}{v}\psi_{1y}^2
+\sum_{i=2}^3\v^{\f12}\f{\mu(\t)}{v}\psi_{iy}^2+\v^{\f12}\f{\lambda(\t)}{v\t}\z_y^2)dy.
\end{array}
 \label{(4.60)}
\end{equation}
Integrating (\ref{(4.58)}) with respect to $y$ yields
\begin{equation}
E_{3\tau}+\f12K_3\leq
C\delta\v^{\f12}(1+\v^{\f12}\tau)^{-1}\|(\Phi_y,\Psi_y,W_y)\|^2+C\delta\v^{\f12}(1+\v^{\f12}\tau)^{-\f32}+\int
N_2dy. \label{(4.61)}
\end{equation}
Here, we only consider the term
$-\int\psi_1\int\xi_1^2\bar\T_{1y}d\xi dy$ because other terms in
$\int N_2dy$ can be estimated similarly. By (\ref{(4.42)}), one has
\begin{equation}
\begin{array}{ll}
&\di -\int\psi_1\int\xi_1^2\bar\T_{1y}d\xi
dy\\[4mm]
&\di =\int\psi_{1y}\int\xi_1^2\bar\T_{1}d\xi dy\\
&\di \leq \f18 K_3
+C\delta\v^{\f12}(1+\v^{\f12}\tau)^{-\f32}+C\v^{\f12}\sum_{|\a^\prime|=1}\int\int\f{\nu(|\xi|)}{\mb{M}_*}|\partial^{\a^\prime}
\bar\mb{G}|^2d\xi
dy\\
&\quad\di+C(\delta+\gamma)\int\int\f{\nu(|\xi|)}{\mb{M}_*}|\bar\mb{G}_1|^2d\xi
dy.\end{array}\label{(4.62)}
\end{equation}
Combining (\ref{(4.61)})) and (\ref{(4.62)}) yields
\begin{equation}
\begin{array}{l}
\di E_{3\tau}+\f14K_3\leq
C_6\delta\v^{\f12}(1+\v^{\f12}\tau)^{-1}\|(\Phi_y,\Psi_y,W_y)\|^2+C_6\delta\v^{\f12}(1+\v^{\f12}\tau)^{-\f32}\\
\quad\di+C_6\v^{\f12}\sum_{|\a^\prime|=1}\int\int\f{\nu(|\xi|)}{\mb{M}_*}|\partial^{\a^\prime}
\bar\mb{G}|^2d\xi
dy+C_6(\delta+\gamma)\int\int\f{\nu(|\xi|)}{\mb{M}_*}|\bar\mb{G}_1|^2d\xi
dy.
\end{array}
\label{(4.63)}
\end{equation}

We need to estimate $\v^{\f12}\|\p_y\|^2$ which is not contained
in  $K_3$. Following the same way as in estimating
$\v^{\f12}\|\Phi_y\|^2$ in the previous subsection, we firstly rewrite the
equation $(\ref{(4.51)})_2$ as
\begin{equation}
\begin{array}{l}
\quad\di\f{4}{3}\v^{\f12}\f{\mu(\bar{\t})}{\bar{v}}\p_{y\tau}-\psi_{1\tau}-\v^{-\f12}(p-\bar{p})_y\\
\di= -\f{4\v^\f12}{3}(\f{\mu({\bar{\t}})}{\bar{v}})_y\psi_{1y}
-\f{4}{3}[(\f{\mu({\t})}{v}-\f{\mu({\bar{\t}})}{\bar{v}})u_{1y}]_y
+\v^{-\f12}R_{1y}+\int\xi_1^2\bar\T_{1y}d\xi,
\end{array}
\label{(4.64)}
\end{equation}
by using the equation of conservation of the mass
$(\ref{(4.51)})_1$.

Since
$$
-(p-\bar p)_y=\v^{\f12}\f{\bar p}{\bar v}\p_y-\v^{\f12}\f{R}{\bar
v}\z_y+(\f{p}{v}-\f{\bar p}{\bar v})v_y-(\f{R}{v}-\f{R}{\bar
v})\t_y,
$$
and
$$
\p_y\psi_{1\tau}=(\p_y\psi_1)_{\tau}-(\p_{\tau}\psi_1)_y+\psi_{1y}^2,
$$
then by multiplying (\ref{(4.64)}) by $\v^{\f12}\p_y$, we get
\begin{equation}
\begin{array}{l}
\di(\f{2\mu(\bar{\t})}{3\bar{v}}\v\p_{y}^2-\v^{\f12}\p_y\psi_1)_{\tau}
+\v^{\f12}\f{\bar p}{\bar v}\p_y^2=
\v(\f{2\mu(\bar{\t})}{3\bar{v}})_{\tau}\p_{y}^2
+\v^{\f12}\psi_{1y}^2+\v^{\f12}\f{R}{\bar
v}\z_y\p_y\\
\quad\di-(\f{p}{v}-\f{\bar p}{\bar v})v_y\p_y+(\f{R}{v}-\f{R}{\bar
v})\t_y\p_y-\f{4}{3}\v(\f{\mu(\bar\t)}{\bar
v})_y\psi_{1y}\p_y\\
\quad\di -\f{4}{3}\v^{\f12}[(\f{\mu({\t})}{v}-\f{\mu({\bar{\t}})}
{\bar{v}})u_{1y}]_y\p_y+R_{1y}\p_y+\v^{\f12}\int\xi_1^2\bar\T_{1y}d\xi\p_y.
\end{array}
\label{(4.65)}
\end{equation}
Integrating (\ref{(4.65)}) with respect to $y$ and using the Cauchy
inequality yield
\begin{equation}
\begin{array}{l}
\di
(\int\f{2\mu(\bar{\t})}{3\bar{v}}\v\p_y^2-\v^{\f12}\p_y\psi_1dy)_{\tau}+\int\v^{\f12}
\f{\bar{p}}{2\bar{v}}\p_y^2dy\\
\di\leq
C_7K_3+C_7\delta\v^{\f12}(1+\v^{\f12}\tau)^{-1}\|(\Phi_y,\Psi_y,W_y)\|^2+C_7\delta\v^{\f32}(1+\v^{\f12}\tau)^{-\f32}+C_7\gamma\v^{\f32}\|\psi_{1yy}\|^2\\
\di+C_7(\delta+\gamma)\v^{\f12}\int\int\f{\nu(|\xi|)}{\mb{M}_*}|\bar\mb{G}_1|^2d\xi
dy+C_7(\delta+\gamma)\v\sum_{|\a^\prime|=
1}\int\int\f{\nu(|\xi|)}{\mb{M}_*}|\partial^{\a^\prime}\bar\mb{G}|^2d\xi
dy\\
\di+C_7\v^{\f32}\sum_{|\a|=2}\int\int\f{\nu(|\xi|)}{\mb{M}_*}|\partial^\a
\bar\mb{G}|^2d\xi dy.
\end{array}
\label{(4.67)}
\end{equation}
Here we have used
\begin{equation}
\begin{array}{l}
\di \v^{\f12}\int |\int\xi_1^2\bar\T_{1y}d\xi|^2 dy\leq
C\v^{\f32}\sum_{|\a|=2} \int\int\f{\nu(|\xi|)}{\mb{M}_*}|\partial^\a
\bar\mb{G}|^2d\xi
dy\\
\quad\di+C(\delta+\gamma)\v\sum_{|\a^\prime|= 1}
\int\int\f{\nu(|\xi|)}{\mb{M}_*}|\partial^{\a^\prime}\bar\mb{G}|^2d\xi
dy
+C(\delta+\gamma)\v^{\f12}\int\int\f{\nu(|\xi|)}{\mb{M}_*}|\bar\mb{G}_1|^2d\xi
dy\\
\quad\di+C\delta\v^2(\v^{\f12}\|\p_y\|^2+K_3)+
C\delta\v^{\f32}(1+\v^{\f12}\tau)^{-\f52}.
\end{array}
\label{(4.68)}
\end{equation}
To estimate $\v^{\f12}\|(\p_{\tau},\psi_{\tau},\z_{\tau})\|^2$, we need to
use the equation (\ref{(4.4+)}). By multiplying $(\ref{(4.4+)})_1$ by
$\v^{\f12}\p_{\tau}$, $(\ref{(4.4+)})_2$ by $\v^{\f12}\psi_{1\tau}$,
$(\ref{(4.4+)})_3$ by $\v^{\f12}\psi_{i\tau}~(i=2,3)$ and
$(\ref{(4.4+)})_4$ by $\v^{\f12}\z_{\tau}$ respectively, and adding
them together, after integrating with respect to $y$, we have
\begin{equation}
\begin{array}{l}
\di \v^{\f12}\|(\p_{\tau},\psi_{\tau},\z_{\tau})\|^2\leq
C_8(\v^{\f12}\|\p_y\|^2+K_3)+C_8\delta\v^{\f12}(1+\v^{\f12}\tau)^{-1}\|(\Phi_y,\Psi_y,W_y)\|^2\\
\quad\di+C_8\delta\v^{\f12}(1+\v^{\f12}\tau)^{-\f32}
+C_8\v^{\f12}\int\int\f{\nu(|\xi|)}{\mb{M}_*}|\bar\mb{G}_y|^2d\xi
dy.
\end{array}
\label{(4.69)}
\end{equation}
Choose large constants $\bar C_4,\bar C_5>1$ such that
$$
\bar C_4E_3+\bar
C_5\int(\f{2\mu(\bar{\t})}{3\bar{v}}\v\p_y^2-\v^{\f12}\p_y\psi_1)dy\geq
\f{\bar C_4}{2}E_3+\bar
C_5\int\f{\mu(\bar{\t})}{3\bar{v}}\v\p_y^2dy,
$$
and
$$
\begin{array}{l}
\di \f{\bar C_4}{4}K_3-\bar C_4C_7K_3-C_3(C_8+1)K_3\geq \f{\bar
C_4}{8}K_3,\\[3mm]
\di\bar C_5\int\v^{\f12}
\f{\bar{p}}{2\bar{v}}\p_y^2dy-C_3(C_8+1)\v^{\f12}\|\p_y\|^2\geq
\f{C_5}{2}\int\v^{\f12} \f{\bar{p}}{2\bar{v}}\p_y^2dy.
\end{array}
$$
Let
\begin{equation}
E_4=\bar C_4E_3+\bar
C_5\int(\f{2\mu(\bar{\t})}{3\bar{v}}\v\p_y^2-\v^{\f12}\p_y\psi_1)dy+\int\int\f{\bar\mb{G}_1^2}{2\mb{M}_*}d\xi
dy,\label{(4.70)}
\end{equation}
\begin{equation}
K_4=\f{\bar C_4}{8}K_3+\f{\bar C_5}{2}\int\v^{\f12}
\f{\bar{p}}{2\bar{v}}\p_y^2dy+\v^{\f12}\|(\p_\tau,\psi_\tau,\z_\tau)\|^2+\f{\sigma}{4}\v^{-\f12}\int\int\f{\nu(|\xi|)}{2\mb{M}_*}|\bar\mb{G}_1|^2d\xi
dy, \label{(4.71)}
\end{equation}
then from (\ref{(4.45)}), (\ref{(4.63)}), (\ref{(4.67)}) and
(\ref{(4.69)}), we have
\begin{equation}
\begin{array}{l}
E_{4\tau}+K_4\leq C_9\delta\v^{\f12}(1+\v^{\f12}\tau)^{-1}\|(\Phi_y,\Psi_y,W_y)\|^2+C_9\delta\v^{\f12}(1+\v^{\f12}\tau)^{-\f32}+C_9\gamma\v^{\f32}\|\psi_{1yy}\|^2\\
\di\qquad\qquad +C_9\v^{\f12}\sum_{|\a^\prime|=
1}\int\int\f{\nu(|\xi|)}{\mb{M}_*}|\partial^{\a^\prime}\bar\mb{G}|^2d\xi
dy+C_9\v^{\f32}\sum_{|\a|=2}\int\int\f{\nu(|\xi|)}{\mb{M}_*}|\partial^\a
\bar\mb{G}|^2d\xi dy.
\end{array}
\label{(4.72)}
\end{equation}

 Next we derive the estimate on
the higher order derivatives. By multiplying $(\ref{(4.51)})_2$
by $-\v\psi_{1yy}$, $(\ref{(4.51)})_3$ by $-\v\psi_{iyy}~(i=2,3)$,
$(\ref{(4.51)})_4$ by $-\v\z_{yy}$, and adding them together, we
obtain
\begin{equation}
\begin{array}{l}
\di (\sum_{i=1}^3\v\f{\psi_{1y}^2}{2}+\v\f{\z_y^2}{2})_{\tau}
+\f{4}{3}\v^{\f32}\f{\mu(\t)}{v}\psi_{1yy}^2+\sum_{i=2}^3\v^{\f32}\f{\mu(\t)}{v}\psi_{iyy}^2
+\v^{\f32}\f{\lambda(\t)}{v}\z_{yy}^2=\\
\di
-\f{4}{3}\v^{\f32}(\f{\mu(\t)}{v})_y\psi_{1y}\psi_{1yy}-\sum_{i=2}^3\v^{\f32}(\f{\mu(\t)}{v})_y\psi_{iy}\psi_{iyy}
-\v^{\f32}(\f{\lambda(\t)}{v})_y\z_y\z_{yy}\\
\di-\f{4}{3}\v[(\f{\mu(\t)}{v}-\f{\mu(\bar\t)}{\bar v})\bar
u_{1y}]_y\psi_{1yy}-\v[(\f{\lambda(\t)}{v}-\f{\lambda(\bar\t)}{\bar
v})\bar\t_y]_y\z_{yy}+\v^{\f12}(p-\bar
p)_y\psi_{1yy}\\
\di +\v^\f12 R_{1y}\psi_{1yy}+\v^\f12(pu_{1y}-\bar p\bar
u_{1y})\z_{yy}-\v Q_5\z_{yy}+\v\sum_{i=1}^3\psi_{iyy}\int\xi_1\xi_i\bar\T_{1y}d\xi\\
\di
-\v\z_{yy}(\sum_{i=1}^3u_i\int\xi_1\xi_i\bar\T_{1y}d\xi-\f12\int\xi_1|\xi|^2\bar\T_{1y}d\xi),
\end{array}
\label{(4.73)}
\end{equation}
where $Q_5$ is defined in (\ref{(4.52)}).

Integrating (\ref{(4.73)}) with respect to $y$ yields
\begin{equation}
\begin{array}{l}
\di
 (\int\v\sum_{i=1}^3\f{\psi_{iy}^2}{2}+\v\f{\z_y^2}{2}dy)_{\tau}+\int(\f{4}{3}\v^{\f32}\f{\mu(\t)}{v}\psi_{1yy}^2+\sum_{i=2}^3\v^{\f32}\f{\mu(\t)}{v}\psi_{iyy}^2
 +\v^{\f32}\f{\lambda(\t)}{v}\z_{yy}^2)dy\\
\di \leq
C(\v^{\f12}\|\p_y\|^2+K_3)+C\delta\v^{\f12}(1+\v^{\f12}\tau)^{-1}\|(\Phi_y,\Psi_y,W_y)\|^2+C\delta\v^{\f32}(1+\v^{\f12}\tau)^{-\f32}\\
\quad\di +
C(\delta+\gamma)\v^{\f12}\int\int\f{\nu(|\xi|)}{\mb{M}_*}|\bar\mb{G}_1|^2d\xi
dy+C(\delta+\gamma)\v\sum_{|\a^{\prime}|=1}\int\int\f{\nu(|\xi|)}{\mb{M}_*}|\partial^{\a^\prime}\bar\mb{G}|^2d\xi
dy\\
\quad\di
+C\v^{\f32}\sum_{|\a|=2}\int\int\f{\nu(|\xi|)}{\mb{M}_*}|\partial^\a
\bar\mb{G}|^2d\xi dy.
\end{array}
\label{(4.74)}
\end{equation}
Now we get the estimation of $\v^{\f32}\|\p_{yy}\|^2$. By applying
$\partial_y$ to $(\ref{(4.4+)})_2$, we get
\begin{equation}
\psi_{1y\tau}+\v^{-\f12}(p-\bar p)_{yy}
=-\f{4}{3}(\f{\mu(\bar\t)}{\bar v}\bar
u_{1y})_{yy}-\v^{-\f12}R_{1yy}-\int\xi_1^2\bar\mb{G}_{yy}d\xi.
\label{(4.75)}
\end{equation}
Note that
\begin{equation}
(p-\bar{p})_{yy}=-\v^{\f12}\f{p}{v}\p_{yy}+\v^{\f12}\f{R}{v}\z_{yy}-\f1v(p-\bar{p})\bar{v}_{yy}
-\v^{\f12}\f{\p}{v}\bar{p}_{yy}-\f{2v_y}{v}(p-\bar{p})_y-\v^{\f12}\f{2\bar{p}_y}{v}\p_y.
\label{(4.76)}
\end{equation}
Multiplying (\ref{(4.75)}) by $-\v^{\f32}\p_{yy}$ and using
(\ref{(4.76)}) imply
\begin{equation}
\begin{array}{l}
\di
-(\int\v^{\f32}\psi_{1y}\p_{yy}dy)_{\tau}+\int\v^{\f32}\f{p}{2v}\p_{yy}^2dy\leq
C\v^{\f32}\|(\psi_{1yy},\z_{yy})\|^2+C\delta\v^{\f12}(1+\v^{\f12}\tau)^{-\f32}\\
\quad\di+C\delta\v^{\f12}(1+\v^{\f12}\tau)^{-1}\|(\Phi_y,\Psi_y,W_y)\|^2+C(\delta+\gamma)(\v^{\f12}\|\p_y\|^2+K_3)\\
\quad\di+C\v^{\f32}
\sum_{|\a|=2}\int\int\f{\nu(|\xi|)}{\mb{M}_*}|\partial^\a
\bar\mb{G}|^2d\xi dy.
\end{array}
\label{(4.77)}
\end{equation}
To estimate $\v^{\f32}\|(\p_{y\tau},\psi_{y\tau},\z_{y\tau})\|^2$
and $\v^{\f32}\|(\p_{\tau\tau},\psi_{\tau\tau},\z_{\tau\tau})\|^2$,
we use the system (\ref{(4.4+)}) again. By applying $\partial_y$ to
(\ref{(4.4+)}), and multiplying the four equations of (\ref{(4.4+)})
by $\v^{\f32}\p_{y\tau}$, $\v^{\f32}\psi_{1y\tau}$,
$\v^{\f32}\psi_{iy\tau}$ $(i= 2,3)$, $\v^{\f32}\z_{y\tau}$
respectively, then adding them together and integrating with respect
to $y$ give
\begin{equation}
\begin{array}{l}
\di
 \v^{\f32}\|(\p_{y\tau},\psi_{y\tau},\z_{y\tau})\|^2\leq
 C\v^{\f32}\|(\p_{yy},\psi_{yy},\z_{yy})\|^2+C\delta\v^{\f12}(1+\v^{\f12}\tau)^{-\f32}\\
\quad\di+C\delta\v^{\f12}(1+\v^{\f12}\tau)^{-1}\|(\Phi_y,\Psi_y,W_y)\|^2+C(\delta+\gamma)(\v^{\f12}\|\p_y\|^2+K_3)\\
\quad\di+C\v^{\f12}
 \int\int\f{\nu(|\xi|)}{\mb{M}_*}|\bar\mb{G}_y|^2d\xi
 dy+C\v^{\f32}
 \sum_{|\a|=2}\int\int\f{\nu(|\xi|)}{\mb{M}_*}|\partial^\a
 \bar\mb{G}|^2d\xi dy.
\end{array}
\label{(4.78)}
\end{equation}
Similarly, we have
\begin{equation}
\begin{array}{l}
\di
 \v^{\f32}\|(\p_{\tau\tau},\psi_{\tau\tau},\z_{\tau\tau})\|^2\leq
 C\v^{\f32}\|(\p_{y\tau},\psi_{y\tau},\z_{y\tau})\|^2+C\delta\v^{\f12}(1+\v^{\f12}\tau)^{-\f32}\\
\quad\di+C\delta\v^{\f12}(1+\v^{\f12}\tau)^{-1}\|(\Phi_y,\Psi_y,W_y)\|^2
 +C(\delta+\gamma)\v^{\f12}\sum_{|\a^\prime|=1}\|\partial^{\a^\prime}(\p,\psi,\z)\|^2\\
\quad\di+C\v^{\f12}
 \int\int\f{\nu(|\xi|)}{\mb{M}_*}|\bar\mb{G}_y|^2d\xi dy+C\v^{\f32}
 \sum_{|\a|=2}\int\int\f{\nu(|\xi|)}{\mb{M}_*}|\partial^\a
 \bar\mb{G}|^2d\xi dy.
\end{array}
\label{(4.79)}
\end{equation}
By choosing $\bar{C}_6$ and $\bar{C}_7$ to be large enough, we have
\begin{equation}
\begin{array}{l}
\di(\bar{C}_6\int\v\sum_{i=1}^3\f{\psi_{iy}^2}{2}+\v\f{\z_y^2}{2}dy-\bar{C}_7\int\v^{\f32}\psi_{1y}\p_{yy}dy)_{\tau}
+\v^{\f32}\sum_{|\a|=2}\|\partial^\a(\p,\psi,\z)\|^2\\
\di \leq
C\v^{\f32}\sum_{|\a|=2}\int\int\f{\nu(|\xi|)}{\mb{M}_*}|\partial^\a
\bar\mb{G}|^2d\xi dy+C\v^{\f12}
\sum_{|\a^\prime|=1}\int\int\f{\nu(|\xi|)}{\mb{M}_*}|\partial^{\a^\prime}
\bar\mb{G}|^2d\xi
dy\\
\di\quad
+C(\delta+\gamma)\v^{\f12}\sum_{|\a^\prime|=1}\|\partial^{\a^\prime}(\p,\psi,\z)\|^2
+C\delta\v^{\f12}(1+\v^{\f12}\tau)^{-1}\|(\Phi_y,\Psi_y,W_y)\|^2\\
\di\quad +C\delta\v^{\f12}(1+\v^{\f12}\tau)^{-\f32}.
\end{array}
\label{(4.80)}
\end{equation}
To close the a priori estimate, we also need to estimate the
derivatives on the non-fluid component $\bar\mb{G}$, i.e.,
$\partial^\a \bar\mb{G}, (|\a|=1,2)$. Applying $\partial_y$ on
(\ref{(3.15)}), we have
\begin{equation}
\begin{array}{l}
\quad\di
\bar\mb{G}_{y\tau}-(\f{u_1}{v}\bar\mb{G}_y)_y+\v^{-\f12}\{\f1v\mb{P}_1(\xi_1\mb{M}_y)\}_y+\{\f1v\mb{P}_1(\xi_1\bar\mb{G}_y)\}_y\\
\di
=\v^{-\f12}\mb{L}_\mb{M}\bar\mb{G}_y+2\v^{-\f12}Q(\mb{M}_y,\bar\mb{G})+2Q(\bar\mb{G}_y,\bar\mb{G}).
\end{array}
\label{(4.81)}
\end{equation}
Since
$$
\mb{P}_1(\xi_1\mb{M}_y)=\f1{Rv\t}\mb{P}_1[\xi_1(\f{|\xi-u|^2}{2\t}\t_y+\xi\cdot
u_y)\mb{M}],
$$
we have
$$
|\{\f1v\mb{P}_1(\xi_1\mb{M}_y)\}_y|\leq
C(v_y^2+u_y^2+\t_y^2+|\t_{yy}|+|u_{yy}|)|\hat{B}(\xi)|\mb{M},
$$
where $\hat{B}(\xi)$ is a polynomial of $\xi$. This yields that
$$
\begin{array}{ll}
\di\v^{\f12}\int\int|\{\f1v\mb{P}_1(\xi_1\mb{M}_y)\}_y\f{\bar\mb{G}_y}{\mb{M}_*}|d\xi
dy \leq &\di
\f{\sigma}{8}\v^{\f12}\int\int\f{\nu(|\xi|)}{\mb{M}_*}|\bar\mb{G}_y|^2d\xi
dy+C\v^{\f32}\|(\psi_{yy},\z_{yy})\|^2\\
& \di
+C(\delta+\gamma)(\v^{\f12}\|\p_y\|^2+K_3)+C\delta\v^{\f12}(1+\v^{\f12}\tau)^{-\f32}.
\end{array}
$$
Thus, multiplying (\ref{(4.81)}) by $\v\f{\bar\mb{G}_y}{\mb{M}_*}$
and using the Cauchy inequality and Lemmas 4.1-4.4, we get
\begin{equation}
\begin{array}{l}
\di(\int\int\v\f{\bar\mb{G}_y^2}{2\mb{M}_*}d\xi
dy)_{\tau}+\f{\sigma}{2}\v^{\f12}\int\int\f{\nu(|\xi|)}{\mb{M}_*}|\bar\mb{G}_y|^2d\xi
dy\leq
C\v^{\f32}\int\int\f{\nu(|\xi|)}{\mb{M}_*}|\bar\mb{G}_{yy}|^2d\xi
dy\\
\quad\di
+C(\delta+\gamma)\int\int\f{\nu(|\xi|)}{\mb{M}_*}|\bar\mb{G}_1|^2d\xi
dy+C(\delta+\gamma)(\v^{\f12}\|\p_y\|^2+K_3)\\
\quad\di
+C\v^{\f32}\|(\p_{yy},\z_{yy})\|^2+C\delta\v^{\f12}(1+\v^{\f12}\tau)^{-\f32}.
\end{array}
\label{(4.82)}
\end{equation}
Similarly,
\begin{equation}
\begin{array}{l}
\di(\int\int\v\f{\bar\mb{G}_{\tau}^2}{2\mb{M}_*}d\xi
dy)_{\tau}+\f{\sigma}{2}\v^{\f12}\int\int\f{\nu(|\xi|)}{\mb{M}_*}|\bar\mb{G}_{\tau}|^2d\xi
dy\leq
C\v^{\f32}\int\int\f{\nu(|\xi|)}{\mb{M}_*}|\bar\mb{G}_{y\tau}|^2d\xi
dy\\
\di+C(\delta+\gamma)\int\int\f{\nu(|\xi|)}{\mb{M}_*}|\bar\mb{G}_1|^2d\xi
dy+C(\delta+\gamma)\v^{\f12}\int\int\f{\nu(|\xi|)}{\mb{M}_*}|\bar\mb{G}_{y}|^2d\xi
dy\\
\di+C\delta\v^{\f12}(1+\v^{\f12}\tau)^{-\f32}+C(\delta+\gamma)\v^{\f12}\sum_{|\a^\prime|=1}\|\partial^{\a^\prime}(\p,\psi,\z)\|^2
+C\v^{\f32}\|(\psi_{y\tau},\z_{y\tau})\|^2.
\end{array}
\label{(4.83)}
\end{equation}
 Finally, we estimate the estimate on the highest
order derivatives, that is, $\int\v^{\f32}\psi_{1y}\p_{yy}dy$ and
$\v^{\f32}\int\int \f{\nu(|\xi|)|\partial^\a
\bar\mb{G}|^2}{\mb{M}_*}d\xi dy$ with $|\a|=2$ in (\ref{(4.80)}). To
do so, it is sufficient to study the estimate for $\v\int\int
\f{|\partial^\a f|^2}{\mb{M}_*}d\xi dy~(|\a|=2)$ because of
(\ref{(4.9)})- (\ref{(4.12)}). For this, from (\ref{(3.16)}) we have
$$
vf_\tau-u_1f_y+\x_1f_y=\v^{-\f12}vQ(f,f)=
v[\mb{L}_\mb{M}\bar\mb{G}+\v^{\f12}Q(\bar\mb{G},\bar\mb{G})].
$$
Applying $\partial^\a$ $(|\a|=2)$ to the above equation gives
\begin{equation}
\begin{array}{ll}
\di v(\partial^\a
f)_{\tau}-v\mb{L}_\mb{M}\partial^\a\bar\mb{G}-u_1(\partial^\a
f)_y+\x_1(\partial^\a
f)_y\\[3mm]
\di =-\partial^\a v f_\tau+\partial^\a u_1
f_y-\sum_{|\a^\prime|=1}[\partial^{\a-\a^\prime}v\partial^{\a^\prime}f_\tau-\partial^{\a-\a^\prime}u_1\partial^{\a^\prime}f_y]\\
\di \quad
+[\partial^\a(v\mb{L}_\mb{M}\bar\mb{G})-v\mb{L}_\mb{M}\partial^\a\bar\mb{G}]+\v^{\f12}\partial^\a[vQ(\bar\mb{G},\bar\mb{G})].
\end{array}
\label{(4.84)}
\end{equation}
Multiplying (\ref{(4.84)}) by $\v\f{\partial^\a
f}{\mb{M}_*}=\v\f{\partial^\a
\mb{M}}{\mb{M}_*}+\v^{\f32}\f{\partial^\a \bar\mb{G}}{\mb{M}_*}$ yields
\begin{equation}
\begin{array}{l}
\quad\di (\v \f{v|\partial^\a f|^2}{2\mb{M}_*})_{\tau}-\v^{\f32}
v\mb{L}_\mb{M}\partial^\a \bar\mb{G}\cdot \f{\partial^\a
\bar\mb{G}}{\mb{M}_*}\\ \di =\v \f{\partial^\a f}{\mb{M}_*}\bigg\{
-\partial^\a v f_\tau+\partial^\a u_1
f_y-\sum_{|\a^\prime|=1}[\partial^{\a-\a^\prime}v\partial^{\a^\prime}f_\tau-\partial^{\a-\a^\prime}u_1\partial^{\a^\prime}f_y]\\
\di \quad
+[\partial^\a(v\mb{L}_\mb{M}\bar\mb{G})-v\mb{L}_\mb{M}\partial^\a\bar\mb{G}]+\v^{\f12}\partial^\a[vQ(\bar\mb{G},\bar\mb{G})]\bigg\}+\v
v\mb{L}_\mb{M}\partial^\a \bar\mb{G}\cdot\f{\partial^\a
\mb{M}}{\mb{M}_*}+(\cdots)_y.
\end{array}
\label{(4.85)}
\end{equation}
We can compute that
$$
\begin{array}{l}
\di \v\int\int |\partial^\a v f_\tau\f{\partial^\a f}{\mb{M}_*}|d\x
dy \\ \di \leq \v\int |\partial^\a
v|\int(|\mb{M}_\tau|+\v^{\f12}|\bar\mb{G}_\tau|)\f{|\partial^\a\mb{M}|+\v^{\f12}|\partial^\a\bar\mb{G}|}{\mb{M}_*}d\x
dy\\
\di \leq
C(\d+\g)\v^{\f32}\|\partial^{\a}(\p,\psi,\z)\|^2+\f{\sigma}{16}\v^{\f32}\int\int
\f{v|\partial^\a\bar\mb{G}|^2}{\mb{M}_*} d\x dy\\
\di \quad +C(\d+\g)^2\v^{\f12}\int\int
\f{|\bar\mb{G}_\tau|^2}{\mb{M}_*} d\x
dy\\
\di
+C(\delta+\gamma)\v\sum_{|\a^\prime|=1}\|\partial^{\a^\prime}(\p,\psi,\z)\|^2+C\delta\v^{\f12}(1+\v^{\f12}\tau)^{-\f32},
\end{array}
$$
and
$$
\begin{array}{l}
\di \v\sum_{|\a^\prime|=1}\int\int
|\partial^{\a-\a^{\prime}}v\partial^{\a^\prime}f_\tau \f{\partial^\a
f}{\mb{M}_*}| d\x dy\\
\di \leq \v\sum_{|\a^\prime|=1}\int|\partial^{\a-\a^{\prime}}v|\int
(|\partial^{\a^\prime}\mb{M}_\tau|+|\partial^{\a^\prime}\mb{G}_\tau|)
\f{|\partial^\a
\mb{M}|+\v^{\f12}|\partial^{\a}\mb{G}|}{\mb{M}_*} d\x dy\\
\di \leq \f{\sigma}{16}\v^{\f32}\int\int
\f{v|\partial^\a\bar\mb{G}|^2}{\mb{M}_*} d\x
dy+C(\d+\g)\v^{\f32}\|\partial^{\a}(\p,\psi,\z)\|^2+C\delta\v^{\f12}(1+\v^{\f12}\tau)^{-\f32}.
\end{array}
$$
Similar estimates can be got to the terms $\v\partial^\a u_1
f_y\f{\partial^\a f}{\mb{M}_*}$ and
$\v\sum_{|\a^\prime|=1}\partial^{\a-\a^\prime}u_1\partial^{\a^\prime}f_y\f{\partial^\a
f}{\mb{M}_*}$.

Also, we have
$$
\begin{array}{l}
\di
\partial^\a(v\mb{L}_\mb{M}\bar\mb{G})-v\mb{L}_\mb{M}\partial^\a\bar\mb{G}=(\partial^\a
v) \mb{L}_\mb{M}\bar\mb{G}+2v Q(\partial^\a\mb{M},\bar\mb{G})\\
\di \qquad
+\sum_{|\a^\prime|=1}\bigg\{2vQ(\partial^{\a-\a^\prime}\mb{M},\partial^{\a^\prime}\bar\mb{G})
+\partial^{\a-\a^\prime}v[\mb{L}_\mb{M}\partial^{\a^\prime}\bar\mb{G}+2Q(\partial^{\a^\prime}\mb{M},\bar\mb{G})]\bigg\},
\end{array}
$$
and
$$
\begin{array}{l}
\di
\v^{\f12}\partial^\a[vQ(\bar\mb{G},\bar\mb{G})]=\v^{\f12}(\partial^\a
v)Q(\bar\mb{G},\bar\mb{G})+\v^{\f12}2vQ(\partial^\a\bar\mb{G},\bar\mb{G})\\
\di\qquad
+\sum_{|\a^\prime|=1}\bigg\{vQ(\partial^{\a-\a^\prime}\bar\mb{G},\partial^{\a^\prime}\bar\mb{G})
+2(\partial^{\a-\a^\prime}v)Q(\partial^{\a^\prime}\bar\mb{G},\bar\mb{G})]\bigg\}.
\end{array}
$$
We only compute one of the above terms as follows, the other terms
can be calculated similarly.
$$
\begin{array}{l}
\quad\di \v^2\int\int\f{v\partial^\a \bar\mb{G}\cdot Q(\partial^\a
\bar\mb{G},\bar\mb{G})}{\mb{M}_*}d\xi dy\\\di
\leq\f{\sigma}{16}\v^{\f32} \int \int\f{v|\partial^\a
\bar\mb{G}|^2}{\mb{M}_*}d\xi dy\\
\qquad\di+C\v^{\f52}\int\bigg(\int\f{\nu(|\xi|)|\partial^\a
\bar\mb{G}|^2}{\mb{M}_*}d\xi \cdot\int\f{|
\bar\mb{G}|^2}{\mb{M}_*}d\xi+\int\f{ |\partial^\a
\bar\mb{G}|^2}{\mb{M}_*}d\xi \cdot\int\f{\nu(|\xi|)|
\bar\mb{G}|^2}{\mb{M}_*}d\xi\bigg) dy\\
\le \di\f{\sigma}{8}\v^{\f32}
\int\int\f{\nu(|\xi|)}{\mb{M}_*}v|\partial^\a \bar\mb{G}|^2d\xi dy
+C\v^{\f52}\sup_{y}\int\f{\nu(|\xi|)|
\bar\mb{G}_1|^2}{\mb{M}_*}d\xi\cdot\int\int \f{|\partial^\a
\bar\mb{G}|^2}{\mb{M}_*}d\xi dy\\
\di\leq  \f{\sigma}{8}\v^{\f32}
\int\int\f{\nu(|\xi|)}{\mb{M}_*}v|\partial^\a \bar\mb{G}|^2d\xi
dy+C(\delta+\g)^2\v^{\f12}\int\int\f{\nu(|\xi|)[|\bar\mb{G}_{1y}|^2+|\bar\mb{G}_1|^2]}{\mb{M}_*}d\xi
dy\\
\di\leq  \f{\sigma}{8}\v^{\f32}
\int\int\f{\nu(|\xi|)}{\mb{M}_*}v|\partial^\a \bar\mb{G}|^2d\xi
dy+C(\delta+\g)^2\v^{\f12}\int\int\f{\nu(|\xi|)[|\bar\mb{G}_{y}|^2+|\bar\mb{G}_1|^2]}{\mb{M}_*}d\xi
dy\\[3mm]
\di\qquad\qquad\qquad
+C\d\v^{\f12}(1+\v^{\f12}\tau)^{-\f32}+C(\d+\g)^2\v^{\f32}\|(\p_y,\psi_y,\z_y)\|^2.
\end{array}
$$

Now we estimate the term $\di \v\int\int v\mb{L}_\mb{M}\partial^\a
\bar\mb{G}\cdot\f{\partial^\a \mb{M}}{\mb{M}_*}d\x dy$ in
(\ref{(4.85)}). Firstly, note that $\mb{P}_1(\partial^\a \mb{M})$
does not contain the term $\partial^\a(v,u,\t)$ for $|\a|=2$. Thus,
we have
\begin{equation}
\begin{array}{l}
\quad\di \v\int\int\f{v\mb{L}_\mb{M}\partial^\a
\bar\mb{G}\cdot\partial^\a \mb{M}}{\mb{M}}d\xi
dy=\v\int\int\f{v\mb{L}_\mb{M}\partial^\a \bar\mb{G}\cdot
\mb{P}_1(\partial^\a \mb{M})}{\mb{M}}d\xi dy\\
\di\leq \f{\sigma }{16}\v^{\f32} \int\int\f{v|\partial^\a
\bar\mb{G}|^2}{\mb{M}_*}d\xi
dy+C(\delta+\gamma)\v\sum_{|\a^\prime|=1}\|\partial^{\a^\prime}(\p,\psi,\z)\|^2+C\delta\v^{\f12}(1+\v^{\f12}\tau)^{-\f32}.
\end{array}
\label{(4.86)}
\end{equation}
Also we can get
\begin{equation}
\begin{array}{l}
\di \v\int\int v\mb{L}_\mb{M}
\partial^\a
\bar\mb{G}\cdot\partial^\a \mb{M}(\f{1}{\mb{M}_*}-\f{1}{\mb{M}})d\xi
dy\leq \f{\sigma }{16}\v^{\f32}
\int\int\f{\nu(|\xi|)}{\mb{M}_*}v|\partial^\a \bar\mb{G}|^2d\xi
dy\\
\quad\di+C\eta_0^2~\v^{\f32}\|\partial^\a(\p,\psi,\z)\|^2+C(\delta+\gamma)
\v^{\f12}\sum_{|\a^\prime|=1}\|\partial^{\a^\prime}(\p,\psi,\z)\|^2+C\delta\v^{\f12}(1+\v^{\f12}\tau)^{-\f32},
\end{array}
\label{(4.87)}
\end{equation}
where the small constant $\eta_0$ is defined in Lemma 4.2. The
combination of (\ref{(4.86)}) and (\ref{(4.87)}) gives the
estimation of $\di \v\int\int v\mb{L}_\mb{M}\partial^\a
\bar\mb{G}\cdot\f{\partial^\a \mb{M}}{\mb{M}_*}d\x dy$.

Thus integrating (\ref{(4.85)}) and recalling all the above
estimates imply
\begin{equation}
\begin{array}{l}
\di (\int\int\v \f{v|\partial^\a f|^2}{2\mb{M}_*}d\xi
dy)_{\tau}+\f{\sigma}{2}\v^{\f32}\int\int\f{\nu(|\xi|)}{\mb{M}_*}v|\partial^\a
\bar\mb{G}|^2d\xi dy\\
\di\leq
C(\delta+\gamma)\v^{\f12}\sum_{|\a^\prime|=1}\|\partial^{\a^\prime}(\p,\psi,\z)\|^2+
C(\eta_0+\delta+\gamma)\v^{\f32}\sum_{|\a|=2}\|\partial^\a(\p,\psi,\z)\|^2\\
\quad\di+C(\delta+\gamma)\v\sum_{|\a^\prime|=1}\int\int\f{\nu(|\xi|)}{\mb{M}_*}|\partial^{\a^\prime}
\bar\mb{G}|^2d\xi
dy+C\delta\v^{\f12}(1+\v^{\f12}\tau)^{-\f32}\\
\quad\di+C(\delta+\gamma)\v^{\f12}\int\int\f{\nu(|\xi|)}{\mb{M}_*}|\bar\mb{G}_1|^2d\xi
dy.
\end{array}
\label{(4.88)}
\end{equation}
By (\ref{(4.9)})-(\ref{(4.11)}), we can choose suitable constants
$\hat{C}_i>1$, $i=1,2,3,4$ so that
\begin{equation}
\begin{array}{l}
\di
E_5=\hat{C}_1E_4+\hat{C}_2(\bar{C}_6\int\v\sum_{i=1}^3\f{\psi_{iy}^2}{2}+\v\f{\z_y^2}{2}dy-\bar
C_7\int\v^{\f32}\psi_{1y}\p_{yy}dy)\\
\qquad\di+\hat{C}_3
\v\sum_{|\a^\prime|=1}\int\int\f{|\partial^{\a^\prime}
\bar\mb{G}|^2}{2\mb{M}_*}d\xi dy+\hat{C}_4
\v\sum_{|\a|=2}\int\int\f{v|\partial^\a f|^2}{2\mb{M}_*}d\xi dy\\
\di\geq C\bigg[
\|(\p,\psi,\z)\|^2+\v\|(\p_y,\psi_y,\z_y)\|^2+\int\int\f{|\bar\mb{G}_1|^2}{\mb{M}_*}d\xi
dy +\v\sum_{|\a^\prime|=1}\int\int\f{|\partial^{\a^\prime}
\bar\mb{G}|^2}{\mb{M}_*}d\xi dy\\
\qquad\di+ \v\sum_{|\a|=2}\int\int\f{|\partial^\a
f|^2}{\mb{M}_*}d\xi dy\bigg]-C\delta\v(1+\v^{\f12}\tau)^{-\f32}.
\end{array}
\label{(4.89)}
\end{equation}
Let
\begin{equation}
\begin{array}{l}
\di
K_5=C^{-1}\bigg[\v^{\f12}\sum_{|\a^\prime|=1}\|\partial^{\a^\prime}(\p,\psi,\z)\|^2+\v^{\f32}\sum_{|\a|=2}\|\partial^\a(\p,\psi,\z)\|^2+\v^{-\f12}\int\int\f{\nu(|\xi|)}{\mb{M}_*}|\bar\mb{G}_1|^2d\xi
dy\\
\di+\v^{\f12}\sum_{|\a^\prime|=1}\int\int\f{\nu(|\xi|)}{\mb{M}_*}|\partial^{\a^\prime}
\bar\mb{G}|^2d\xi
dy+\v^{\f32}\sum_{|\a|=2}\int\int\f{\nu(|\xi|)}{\mb{M}_*}|\partial^\a
\bar\mb{G}|^2d\xi dy\bigg].
\end{array}
\label{(4.90)}
\end{equation}
Then by the estimates (\ref{(4.72)}), (\ref{(4.80)}),
(\ref{(4.82)}), (\ref{(4.83)}), (\ref{(4.88)}), we obtain
\begin{equation}
E_{5\tau}+K_5\leq
C\delta\v^{\f12}(1+\v^{\f12}\tau)^{-\f32}+C\delta\v^{\f12}(1+\v^{\f12}\tau)^{-1}\|(\Phi_y,\Psi_y,W_y)\|^2.
\label{(4.91)}
\end{equation}

\section{The proof of Theorem 3.1}
\setcounter{equation}{0}
Choose a large constant $\hat{C}_5$ and set
\begin{equation}
E_6=E_2+\hat{C}_5E_5, \qquad K_6=K_2+\hat{C}_5K_5.\label{(5.1)}
\end{equation}
By combining (\ref{(4.49)}) and (\ref{(4.91)}), we have
\begin{equation}
\begin{array}{ll}
&\di E_{6\tau}+K_6\\ & \di\leq
C_0\delta\v^{\f12}(1+\v^{\f12}\tau)^{-1}E_2+C_0\delta\v^{\f12}(1+\v^{\f12}\tau)^{-1}\|(\Phi_y,\Psi_y,W_y)\|^2+C_0\delta\v^{\f12}(1+\v^{\f12}\tau)^{-\f12}\\
 & \di\leq C_0\delta\v^{\f12}(1+\v^{\f12}\tau)^{-1}E_6+C_0\delta\v^{\f12}(1+\v^{\f12}\tau)^{-\f12}.
\end{array}
\label{(5.2)}
\end{equation}
Then Gronwall inequality implies that
\begin{equation}
E_6(\tau)\leq C(E_6(0)+\delta)(1+\v^{\f12}\tau)^{\f12},\quad
\int_0^\tau K_6(y,s)ds\leq C(E_6(0)+\delta)(1+\v^{\f12}\tau)^{\f12}.
\label{(5.3)}
\end{equation}
Now multiplying (\ref{(4.91)}) by $(1+\v^{\f12}\tau)$ gives
\begin{equation}
\begin{array}{ll}
[(1+\v^{\f12}\tau)E_5]_\tau&\di\leq (1+\v^{\f12}\tau)E_{5\tau}+\v^{\f12}E_5\\
&\di \leq
C\delta\v^{\f12}(1+\v^{\f12}\tau)^{-\f12}+C\delta\v^{\f12}\|(\Phi_y,\Psi_y,W_y)\|^2+\v^{\f12}E_5\\
&\di \leq C\delta\v^{\f12}(1+\v^{\f12}\tau)^{-\f12}+C K_6.
\end{array}
\label{(5.4)}
\end{equation}
Integrating (\ref{(5.4)}) with respect to $\tau$ and using
(\ref{(5.3)}) yield that
$$
E_5(\tau)\leq C(E_6(0)+\delta)(1+\v^{\f12}\tau)^{-\f12}.
$$
Thus,  we have
$$
\|(\Phi,\Psi,W)\|_{L^\i_y}^2\leq
C\|(\Phi,\Psi,W)\|\|(\Phi_y,\Psi_y,W_y)\|\leq
CE_6^{\f12}E_5^{\f12}\leq C(E_6(0)+\delta),
$$
and
\begin{equation}
\begin{array}{l}
\di
\|(\p,\psi,\z)\|^2+\v\|(\p_y,\psi_y,\z_y)\|^2+\int\int\f{\bar\mb{G}_1^2}{\mb{M}_*}d\xi
dy\\
\qquad\di+\v\sum_{|\a^\prime|=1} \int\int\f{|\partial^{\a^\prime}
\bar\mb{G}|^2}{\mb{M}_*}d\xi
dy+\v\sum_{|\a|=2}\int\int\f{|\partial^\a
f|^2}{\mb{M}_*}d\xi dy\\
\qquad\di\leq
C(E_6(0)+\delta)(1+\v^{\f12}\tau)^{-\f12}\\
\qquad \di \leq C(E_6(0)+\delta).
\end{array}
\label{(5.5)}
\end{equation}
And this  closes the a priori estimate (\ref{(4.1)}).

Now it remains to prove  the decay rate of (\ref{(3.18)}). By
(\ref{(5.5)}), we have
$$
\v^{\f12}\|(\p,\psi,\z)\|^2_{L^\i_y}\leq
C\v^{\f12}\|(\p,\psi,\z)\|\|(\p_y,\psi_y,\z_y)\|\leq
C(E_6(0)+\delta),
$$
and
$$
\begin{array}{ll}
\di \v^{\f12}\|\int\f{\bar\mb{G}^2}{\mb{M}_*}d\xi\|_{L_\i^y}&\di
\leq C\v^{\f12}\bigg(\int\int\f{\bar\mb{G}^2}{\mb{M}_*}d\xi
dy\bigg)^{\f{1}{2}}\cdot\bigg(\int\int\f{|\bar\mb{G}_{y}|^2}{\mb{M}_*}d\xi
dy\bigg)^{\f{1}{2}}\\[3mm]
&\di \leq C\bigg[\bigg(\int\int\f{|\bar\mb{G}_{1}|^2}{\mb{M}_*}d\xi
dy\bigg)^{\f{1}{2}}+\|(\bar\t_{y},\bar u_{y})\|\bigg](E_6(0)+\d)^{\f12}\\
&\di \leq C(E_6(0)+\d).
\end{array}
$$
Finally,
\begin{equation}
\begin{array}{l}
\di
\sup_y\int\f{|f(y,\tau,\xi)-\mb{M}_{[\bar{v},\bar{u},\bar{\t}]}(y,\tau,\xi)|^2}{\mb{M}_*}d\xi
\\
\di\leq
C\sup_y\int\f{|\mb{M}(y,\tau,\xi)-\mb{M}_{[\bar{v},\bar{u},\bar{\t}]}(y,\tau,\xi)|^2}{\mb{M}_*}d\xi
+\sup_y\int\f{\mb{G}^2}{\mb{M}_*}d\xi \\
\di \leq
C\v\|(\p,\psi,\z)\|_{L_\i}^2+C\v\sup_y\int\f{\bar\mb{G}^2}{\mb{M}_*}d\xi\\
\di \leq C(E_6(0)+\d)\v^{\f12},
\end{array}
\label{(5.6)}
\end{equation}
which gives (\ref{(3.18)}). And this completes the proof of Theorem
3.1.

\end{document}